\documentclass[11pt, a4paper]{amsart}
\usepackage{amsthm}
\theoremstyle{}
{\theoremstyle{definition}
\newtheorem{dfn}{Definition}[section]}
\newtheorem{prop}[dfn]{Proposition}

\newtheorem{thm}[dfn]{Theorem}
{\theoremstyle{definition}
\newtheorem{rem}[dfn]{Remark}}
\newtheorem{lem}[dfn]{Lemma}

\usepackage[top=30truemm,bottom=30truemm,left=30truemm,right=30truemm]{geometry}
\usepackage{amsmath,amssymb,amscd,bm,graphicx}
\usepackage[all]{xy}
\begin{document}

\title[]{Derived Kn\"orrer Periodicity and Orlov's theorem for gauged Landau-Ginzburg models}

\author[Y.~Hirano]{Yuki Hirano}
\date{}


\begin{abstract} We prove a Kn\"orrer periodicity type equivalence between derived factorization categories of gauged LG models, which is an analogy of a theorem proved by Shipman and Isik independently. 
As an application, we obtain a gauged LG version of  Orlov's theorem describing a relationship between categories of graded matrix factorizations and derived categories of hypersurfaces in projective spaces, by combining the above Kn\"orrer periodicity type equivalence and  the theory of  variations of GIT quotients due to Ballard, Favero and Katzarkov. 

\end{abstract}

\maketitle{}
\section{Introduction}

\subsection{Backgrounds and motivations}
When $X$ is a scheme, $G$ is an affine algebraic group acting on $X$,  $\chi:G\rightarrow \mathbb{G}_m$ is a character, and   $W:X\rightarrow \mathbb{A}^1$ is a $\chi$-semi invariant regular function, we call data  $(X,\chi,W)^G$   a {\it gauged Landau-Ginzburg (LG) model}.  Following Positselski \cite{posi}, \cite{efi-posi}, we consider  {\it the derived factorization category} of $(X,\chi,W)^G$, denoted by
$${\rm Dcoh}_G(X,\chi,W).$$
Derived factorization categories are simultaneous generalizations of bounded derived categories of coherent sheaves on schemes, and of categories of (graded) matrix factorizations of (homogeneous) polynomials. 
Orlov proved the following semi-orthogonal decompositions between bounded derived categories of hypersurfaces in projective spaces and categories of graded matrix factorizations \cite{orlov3}. 
\begin{thm}[\cite{orlov3} Theorem 40]
Let $X\subset \mathbb{P}^{N-1}_k$ be the hypersurface defined by a section $f\in \Gamma(\mathbb{P}^{N-1}_k,\mathcal{O}(d))$. Denote by $F$ the corresponding homogeneous polynomial.\\
$(1)$ If $d<N$, there is a semi-orthogonal decomposition 
$${\rm D^b}({\rm coh}X)=\langle \mathcal{O}_X(d-N+1),...,\mathcal{O}_X,{\rm Dcoh}_{\mathbb{G}_m}(\mathbb{A}^N_k,\chi_{d},F)\rangle.$$
$(2)$ If $d=N$, there is an equivalence
$${\rm D^b}({\rm coh}X)\cong {\rm Dcoh}_{\mathbb{G}_m}(\mathbb{A}^N_k,\chi_{d},F).$$
$(3)$ If $d>N$, there is a semi-orthogonal decomposition
$${\rm Dcoh}_{\mathbb{G}_m}(\mathbb{A}^N_k,\chi_{d},F)=\langle k,...,k(N-d+1),{\rm D^b}({\rm coh}X)\rangle.$$
\end{thm}
While Orlov's approach  was algebraic, there are geometric approaches to the above  theorem \cite{shipman}, \cite{vgit}, \cite{degree d}, where a version of Kn\"orrer  periodicity \cite{shipman}, \cite{isik} and homological variations of GIT quotients \cite{segal}, \cite{h-l}, \cite{vgit} are the main tools. Combinations of Kn\"orrer periodicity and the theory of variations of GIT quotients also imply homological projective dualities \cite{hpd via vgit}, \cite{ads}, \cite{seg-tho}, \cite{ren}.

In this article, we prove another version of Kn\"orrer periodicity \cite{knorrer}, which is a derived (or global) version, and we combine it with the theory of variations of GIT quotients by \cite{vgit} to obtain a gauged LG version of Orlov's theorem.
 
 \subsection{Main results}
Let $k$ be an algebraically closed field  of characteristic zero.
 Let $X$ be a smooth quasi-projective variety over $k$, and let $G$ be a reductive affine algebraic group acting on $X$. Let $\mathcal{E}$ be a $G$-equivariant locally free sheaf of finite rank, and choose a $G$-invariant regular section $s\in\Gamma(X,\mathcal{E}^{\vee})^G$. Denote by $Z\subset X$ the zero scheme of $s$. Let $\chi:G\rightarrow \mathbb{G}_m$ be a character of $G$, and set $\mathcal{E}(\chi):=\mathcal{E}\otimes\mathcal{O}(\chi),$ where $\mathcal{O}(\chi)$ is the $G$-equivariant invertible sheaf corresponding to $\chi$.
 Then $\mathcal{E}(\chi)$ induces a vector bundle  ${\rm V}(\mathcal{E}(\chi))$ over $X$ with a $G$-action induced by the equivariant structure of $\mathcal{E}(\chi)$. Let $q:{\rm V}(\mathcal{E}(\chi))\rightarrow X$ and $p:{\rm V}(\mathcal{E}(\chi))|_Z\rightarrow Z$ be natural projections, and let $i:{\rm V}(\mathcal{E}(\chi))|_Z\rightarrow{\rm V}(\mathcal{E}(\chi))$ be a natural inclusion.
 The regular section $s$ induces a $\chi$-semi invariant regular function $Q_s:{\rm V}(\mathcal{E}(\chi))\rightarrow \mathbb{A}^1$. The first main result in this paper is the following:
 \begin{thm}
 Let $W:X\rightarrow \mathbb{A}^1$ be a $\chi$-semi invariant regular function, such that the restricted function $W|_Z:Z\rightarrow \mathbb{A}^1$ is flat. Then there is an equivalence
 $$i_*p^*:{\rm Dcoh}_G(Z,\chi,W|_Z)\xrightarrow{\sim}{\rm Dcoh}_G({\rm V}(\mathcal{E}(\chi)),\chi, q^*W+Q_s).$$
 \end{thm}

The above result  is an analogy of  Shipman's and Isik's result, where they consider the case when $G=\mathbb{G}_m$, the $G$-action on $X$ is trivial, $\chi={\rm id}_{\mathbb{G}_m}$, and $W=0$ \cite{shipman}, \cite{isik}. Furthermore, the above theorem can be considered as a generalization of  Kn\"orrer periodicity  to a derived and $G$-equivariant version. The proof of the above theorem is quite different from Shipman's and Isik's proofs, and we consider relative singularity categories introduced in \cite{efi-posi}, which are equivalent to derived factorization categories, and  use results in \cite{orlov1}. 
\vspace{2mm}

To state the next result, let $S$ be a smooth quasi-projective variety over $k$ with a $\mathbb{G}_m$-action, and let $W:S\rightarrow \mathbb{A}^1$ be a $\chi_1:={\rm id}_{\mathbb{G}_m}$-semi invariant regular function which is flat. 
Let $d>1$ and $N>0$ be positive integers, and consider $\mathbb{G}_m$-actions on $\mathbb{A}^N_S:=S\times \mathbb{A}^N_k$ and on $\mathbb{P}^{N-1}_S:=S\times \mathbb{P}^{N-1}_k$ given by 
$$\mathbb{G}_m\times\mathbb{A}^N_S\ni t\times (s,v_1,...,v_N)\mapsto (t^d\cdot s, tv_1,...tv_N)\in\mathbb{A}_S^N$$
$$\mathbb{G}_m\times\mathbb{P}^{N-1}_S\ni t\times (s,v_1:...:v_N)\mapsto (t\cdot s, v_1:...:v_N)\in\mathbb{P}^{N-1}_S.$$ 
Denote by the same notation $W:\mathbb{A}^N_S\rightarrow\mathbb{A}^1$ and $W:\mathbb{P}^{N-1}_S\rightarrow\mathbb{A}^1$ the pull-backs of $W:S\rightarrow \mathbb{A}^1$ by the natural projections respectively.
  Combining  the above derived Kn\"orrer periodicity with the theory of variations of GIT quotients,  we obtain the following gauged LG version of the Orlov's theorem:

 \begin{thm}
 Let $X\subset\mathbb{P}_S^{N-1}$ be the hypersurface defined by  a  $\mathbb{G}_m$-invariant section $f\in \Gamma(\mathbb{P}_S^{N-1},\mathcal{O}(d))^{\mathbb{G}_m}$, and assume that the morphism $W:\mathbb{P}^{N-1}_S\rightarrow\mathbb{A}^1$ is flat on $X$.
 Denote by $F:\mathbb{A}^N_S\rightarrow \mathbb{A}^1$ the regular function induced by $f$. \\
 $(1)$ If $d<N$, there are fully faithful functors 
 $$\Phi:{\rm Dcoh}_{\mathbb{G}_m}(\mathbb{A}^N_S,\chi_d,W+F)\rightarrow {\rm Dcoh}_{\mathbb{G}_m}(X,\chi_1,W)$$
 $$\Upsilon:{\rm Dcoh}_{\mathbb{G}_m}(S,\chi_1,W)\rightarrow {\rm Dcoh}_{\mathbb{G}_m}(X,\chi_1,W),$$
 and there is a semi-orthogonal decomposition
$${\rm Dcoh}_{\mathbb{G}_m}(X,\chi_1,W)=\langle\Upsilon_{d-N+1},...,\Upsilon_0,\Phi( {\rm Dcoh}_{\mathbb{G}_m}(\mathbb{A}^N_S,\chi_d,W+F))\rangle,$$
where $\Upsilon_i$ denotes the essential image of the composition $(-)\otimes\mathcal{O}(i)\circ\Upsilon$.\\

\begin{flushleft}$(2)$ If $d=N$, there is an equivalence\end{flushleft}
$${\rm Dcoh}_{\mathbb{G}_m}(X,\chi_1,W)\cong{\rm Dcoh}_{\mathbb{G}_m}(\mathbb{A}^N_S,\chi_d,W+F).$$\\
$(3)$ If $d>N$, there are fully faithful functors 
 $$\Psi:{\rm Dcoh}_{\mathbb{G}_m}(X,\chi_1,W)\rightarrow {\rm Dcoh}_{\mathbb{G}_m}(\mathbb{A}^N_S,\chi_d,W+F)$$
 $$\Upsilon:{\rm Dcoh}_{\mathbb{G}_m}(S,\chi_1,W)\rightarrow {\rm Dcoh}_{\mathbb{G}_m}(\mathbb{A}^N_S,\chi_d,W+F),$$ and there is a semi-orthogonal decomposition
$${\rm Dcoh}_{\mathbb{G}_m}(\mathbb{A}^N_S,\chi_d,W+F)=\langle\Upsilon_0,...,\Upsilon_{N-d+1},\Psi({\rm Dcoh}_{\mathbb{G}_m}(X,\chi_1,W))\rangle,$$
 where $\Upsilon_i$ denotes the essential image of the composition $(-)\otimes\mathcal{O}(\chi_i)\circ\Upsilon$.
 
  \end{thm} 
 
 \vspace{3mm}
Since we have an equivalence $${\rm D^b}({\rm coh}X)\cong{\rm Dcoh}_{\mathbb{G}_m}(X,\chi_1,0),$$ where the $\mathbb{G}_m$-action on $X$ is trivial, we can view Orlov's theorem as the case when  $S={\rm Spec}k$ and $W=0$ in the above theorem.

\subsection{Notation and conventions}
\begin{itemize}
\item{} For an integer $n\in\mathbb{Z}$, we denote by $\chi_n:\mathbb{G}_m\rightarrow \mathbb{G}_m$ the character of the 1-dimensional algebraic torus $\mathbb{G}_m$ defined by 
$$\chi_n(t):=t^n.$$
\item{} For a character $\chi:G\rightarrow \mathbb{G}_m$ of an algebraic group $G$, we denote by $\mathcal{O}(\chi)$ the $G$-equivariant invertible  sheaf  $(\mathcal{O}_X,\lambda:\pi^*\mathcal{O}_X\xrightarrow{\sim}\sigma^*\mathcal{O}_X)$ associated to $\chi$, where $\pi:G\times X\to X$ and $\sigma:G\times X\to X$ are the projection and the  morphism defining the $G$-action respectively. For any $g\in G$, $\lambda_g:=\lambda|_{\{g\}\times X}:\mathcal{O}_X\xrightarrow{\sim}g^*\mathcal{O}_X$ is  given as the composition $$\mathcal{O}_X\xrightarrow{\chi(g)}\mathcal{O}_X\xrightarrow{\sim}g^*\mathcal{O}_X$$ of  the multiplication by $\chi(g)\in\mathbb{G}_m$ and  the natural isomorphism $\mathcal{O}_X\xrightarrow{\sim}g^*\mathcal{O}_X$.
For a $G$-equivariant  quasi-coherent sheaf $F$ on a $G$-scheme, we set $$F(\chi):=F\otimes \mathcal{O}(\chi).$$
\item{} Throughout this article, unless stated otherwise, all schemes and categories are over an algebraically closed field $k$ of characteristic zero.

\end{itemize}

\subsection{Acknowledgments}
The main part of this article was written whilst the author was visiting the School of Mathematics at  the University of Edinburgh during the period from October 2015 to March 2016. He appreciates the hospitality, and he would like to express his deep gratitude to Michael Wemyss for fruitful discussions and many suggestions to improve  this paper. The author is also grateful to his supervisor Hokuto Uehara, Timothy Logvinenko and Ed Segal for their valuable advices. Finally, I would like to thank an anonymous referee for valuable suggestions for the improvements of the paper.
The author is a Research Fellow of Japan Society for the Promotion of Science. He is partially supported by Grant-in-Aid for JSPS Fellows $\sharp$26-6240.

\section{Derived factorization categories}
In this section, we provide the definition of  the derived factorization category, which is introduced by Positselski (cf. \cite{posi}, \cite{efi-posi}), and its basic properties.
\subsection{Derived factorization categories}\hfill

\begin{dfn}
A \textbf{gauged Landau-Ginzburg model}, or \textbf{ gauged LG model}, is data $(X,\chi,W)^G$ with
$X$ a scheme, $\chi:G\rightarrow\mathbb{G}_m$ a character of $G$, $G$ affine algebraic group acting on $X$, and $W:X\rightarrow\mathbb{A}^1$ $\chi$-semi invariant regular function, i.e.  $W(g\cdot x)=\chi(g)W(x)$ for any $g\in G$ and any $x\in X$. If $G$ is trivial, we denote the gauged LG model by $(X,W)$, and call it \textbf{Landau-Ginzburg model} or \textbf{LG model}. 

\end{dfn}
 
 For a gauged LG model, we consider its factorizations which are $^{\rotatebox[origin=C]{180}{"}}$twisted"  complexes.

\begin{dfn}
Let $(X,\chi,W)^G$ be a gauged LG model. A \textbf{factorization} $F$ of $(X,\chi,W)^G$ is a sequence
$$F=\Bigl(F_1\xrightarrow{\varphi_1^F} F_0\xrightarrow{\varphi_0^F} F_1(\chi)\Bigr),$$
where $F_i$ is a $G$-equivariant quasi-coherent sheaf on $X$ and $\varphi_i^F$ is a $G$-invariant homomorphism for $i=0,1$ such that $\varphi_0^F\circ\varphi_1^F=W\cdot {\rm id}_{F_1}$ and $\varphi_1^F(\chi)\circ\varphi_0^F=W\cdot {\rm id}_{F_0}$. Equivariant quasi-coherent sheaves $F_0$ and $F_1$ in the above sequence are called \textbf{components} of the factorization $F$.
\end{dfn}

\begin{dfn}
For a gauged LG model $(X,\chi,W)^G$, we define a dg category
$${\rm Qcoh}_G(X,\chi,W)$$
whose objects are factorizations of $(X,\chi,W)^G$, and whose complexes of morphisms are defined as follows:
For two objects $E,F\in {\rm Qcoh}_G(X,\chi,W)$, we define the complex ${\rm Hom}(E,F)^{\text{\tiny{\textbullet}}}$ of morphisms from $E$ to $F$ as the following  graded vector space
$${\rm Hom}(E,F)^{\text{\tiny{\textbullet}}}:=\bigoplus_{n\in\mathbb{Z}}{\rm Hom}(E,F)^n$$
with a differential $d^i:{\rm Hom}(E,F)^i\rightarrow{\rm Hom}(E,F)^{i+1}$ given by
$$d^i(f):=\varphi^F\circ f-(-1)^i f\circ\varphi^E$$
where 
$${\rm Hom}(E,F)^{2m}:={\rm Hom}(E_1,F_1(\chi^m))\oplus{\rm Hom}(E_0,F_0(\chi^m))$$
$${\rm Hom}(E,F)^{2m+1}:={\rm Hom}(E_1,F_0(\chi^m))\oplus{\rm Hom}(E_0,F_1(\chi^{m+1})).$$

We define dg full subcategories ${\rm coh}_G(X,\chi,W)$, ${\rm Inj}_G(X,\chi,W)$, ${\rm LFr}_G(X,\chi,W)$, and ${\rm lfr}_G(X,\chi,W)$
 of  ${\rm Qcoh}_G(X,\chi,W)$ whose objects are factorizations whose components are coherent, injective quasi-coherent,  locally free, and locally free of finite ranks respectively.
 If $G$ is trivial, dropping $G$ and $\chi$ from each notation, we denote the dg categories by ${\rm Qcoh}(X,W)$, ${\rm coh}(X,W)$, etc. 
\end{dfn}

\vspace{3mm}The dg category ${\rm Qcoh}_G(X,\chi,W)$ induces two categories $$Z^0({\rm Qcoh}_G(X,\chi,W))$$ $$H^0({\rm Qcoh}_G(X,\chi,W)).$$ Objects of these categories are same as objects of ${\rm Qcoh}_G(X,\chi,W)$, and sets of morphisms are defined as 
$${\rm Hom}_{Z^0({\rm Qcoh}_G(X,\chi,W))}(E,F):=Z^0({\rm Hom}(E,F)^{\text{\tiny{\textbullet}}})={\rm Ker}(d^0)$$
$${\rm Hom}_{H^0({\rm Qcoh}_G(X,\chi,W))}(E,F):=H^0({\rm Hom}(E,F)^{\text{\tiny{\textbullet}}})={\rm Ker}(d^0)/{\rm Im}(d^{-1}).$$

\begin{rem}
 We can write down the above sets of morphisms  as follows:

${\rm Hom}_{Z^0({\rm Qcoh}_G(X,\chi,W))}(E,F)$ is the set of pairs $(f_1,f_0)$ such that $f_i\in{\rm Hom}_{{\rm Qcoh}_GX}(E_i, F_i)$ and that  the following diagram is commutative
\[\xymatrix{
 E_1\ar[rr]^{\varphi_1^E}\ar[d]_{f_1}&&E_0\ar[rr]^{\varphi_0^E}\ar[d]^{f_0}&&E_1(\chi)\ar[d]^{f_1(\chi)}\\
F_1\ar[rr]^{\varphi_1^F}&&F_0\ar[rr]^{\varphi_0^F}&&F_1(\chi).
 }\]

The set of morphisms in the category $H^0({\rm Qcoh}_G(X,\chi,W))$ can be described as the set of homotopy equivalence (denoted by $\sim$) classes in {\rm Hom}$_{Z^0({\rm Qcoh}_G(X,\chi,W))}(E,F)$,
$${\rm Hom}_{H^0({\rm Qcoh}_G(X,\chi,W))}(E,F)={\rm Hom}_{Z^0({\rm Qcoh}_G(X,\chi,W))}(E,F)/\sim,$$
where two morphisms $f=(f_1,f_0)$ and $g=(g_1,g_0)$ in ${\rm Hom}_{Z^0({\rm Qcoh}_G(X,\chi,W))}(E,F)$ are \textbf{homotopy equivalent} if there exist two morphisms 
$$h_0:E_0\rightarrow F_1\hspace{3mm}{\rm and}\hspace{3mm}h_1:E_1(\chi)\rightarrow F_0$$
such that $f_0=\varphi_1^Fh_0+h_1\varphi_0^E$ and $f_1(\chi)=\varphi_0^Fh_1+h_0(\chi)\varphi_1^E(\chi)$.

\end{rem}

\vspace{3mm}
We easily see the following result (cf. \cite[Proposition 3.5]{hirano}). 
\begin{prop}
The categories $Z^0({\rm Qcoh}_G(X,\chi,W))$ and $Z^0({\rm coh}_G(X,\chi,W))$ are abelian, and  the categories $Z^0({\rm LFr}_G(X,\chi,W))$ and $Z^0({\rm lfr}_G(X,\chi,W))$ are exact.
\end{prop}

We next define the totalizations of bounded complexes of factorizations.

\begin{dfn}
Let $F^{\text{\tiny{\textbullet}}}=(\cdot\cdot\cdot\rightarrow F^i\xrightarrow{\delta^i}F^{i+1}\rightarrow\cdot\cdot\cdot)$ be a bounded complex of $Z^0({\rm Qcoh}_G(X,\chi,W))$. For $l=0,1$, set
$$T_l:=\bigoplus_{i+j=-l}F^i_{\overline{j}}(\chi^{\lceil j/2\rceil}),$$
and let $$t_l:T_l\rightarrow T_{\overline{l+1}}$$
be a $G$-invariant homomorphism given by
$$t_l|_{F^i_{\overline{j}}(\chi^{\lceil j/2\rceil})}:=\delta^i_{\overline{j}}(\chi^{\lceil j/2\rceil})+(-1)^i\varphi^{F^i}_{\overline{j}}(\chi^{\lceil j/2\rceil}),$$
where $\overline{n}$ is $n$ modulo $2$, and $\lceil m\rceil$ is the minimum integer which is greater than or equal to a real number $m$.
We define the {\bf totalization} Tot$(F^{\text{\tiny{\textbullet}}})\in Z^0({\rm Qcoh}_G(X,\chi,W))$ of $F^{\text{\tiny{\textbullet}}}$ as 
$${\rm Tot}(F^{\text{\tiny{\textbullet}}}):=\Bigl(T_1\xrightarrow{t_1}T_0\xrightarrow{t_0}T_1(\chi)\Bigr).$$
\end{dfn}

\vspace{3mm}In what follows, we will recall that the category $H^0({\rm Qcoh}_G(X,\chi,W))$ has a structure of a triangulated category.

\vspace*{2mm}
\begin{dfn}
We define an automorphism $T$ on $H^0({\rm Qcoh}_G(X,\chi,W))$, which is called \textbf{shift functor}, as follows.
For an object $F\in H^0({\rm Qcoh}_G(X,\chi,W))$, we define an object $T(F)$ as
$$T(F):=\Bigl(F_0\xrightarrow{-\varphi^F_0}F_1(\chi)\xrightarrow{-\varphi^F_1(\chi)}F_0(\chi)\Bigr)$$
and for a morphism $f=(f_1,f_0)\in {\rm Hom}(E,F)$, we set $T(f):=(f_0,f_1(\chi))\in{\rm Hom}(T(E),T(F))$. For any integer $n\in\mathbb{Z}$, denote by $(-)[n]$ the functor $T^n(-)$.
\end{dfn}

\vspace*{2mm}
\begin{dfn}
Let $f : E\rightarrow F$ be a morphism in $Z^0({\rm Qcoh}_G(X,\chi,W))$. We define its {\bf mapping cone} Cone$(f)$ to be the totalization of the complex $$(\cdot\cdot\cdot\rightarrow0\rightarrow E\xrightarrow{f} F\rightarrow0\rightarrow\cdot\cdot\cdot)$$ with $F$ in degree zero.

A $\textbf{distinguished triangle}$ is a sequence in $H^0({\rm Qcoh}_G(X,\chi,W))$ which is isomorphic to a sequence of the form

$$E\xrightarrow{f}F\xrightarrow{i}{\rm Cone}(f)\xrightarrow{p}E[1],$$
where $i$ and $p$ are  natural injection and  projection respectively.
\end{dfn}

The following is  well-known.

\vspace*{2mm}
\begin{prop}
$H^0({\rm Qcoh}_G(X,\chi,W))$ is a triangulated category with respect to its shift functor  and its distinguished triangles defined above. Full subcategories $H^0({\rm coh}_G(X,\chi,W))$, $H^0({\rm Inj}_G(X,\chi,W))$, $H^0({\rm LFr}_G(X,\chi,W))$, and $H^0({\rm lfr}_G(X,\chi,W))$ are full triangulated subcategories.
\end{prop}

\vspace*{2mm}
Following Positselski (\cite{posi}, \cite{efi-posi}), we define derived factorization categories.

\begin{dfn}
Denote by  ${\rm Acycl}({\rm coh}_G(X,\chi,W))$ the smallest thick subcategory of \\$H^0({\rm coh}_G(X,\chi,W))$ containing all totalizations of short exact sequences in $Z^0({\rm coh}_G(X,\chi,W))$.  We define the \textbf{derived factorization category} of $(X,\chi,W)^G$ as the Verdier quotient 
$${\rm Dcoh}_G(X,\chi,W):=H^0({\rm coh}_G(X,\chi,W))/{\rm Acycl}({\rm coh}_G(X,\chi,W)).$$
\end{dfn}

\vspace{3mm}
Similarly, consider thick full subcategories  ${\rm Acycl}({\rm Qcoh}_G(X,\chi,W))$, ${\rm Acycl}({\rm LFr}_G(X,\chi,W))$ and ${\rm Acycl}({\rm lfr}_G(X,\chi,W))$ of $H^0({\rm Qcoh}_G(X,\chi,W))$, $ H^0({\rm LFr}_G(X,\chi,W))$, and $H^0({\rm lfr}_G(X,\chi,W))$ respectively, and 
denote the corresponding Verdier quotients by
$${\rm DQcoh}_G(X,\chi,W):=H^0({\rm Qcoh}_G(X,\chi,W))/{\rm Acycl}({\rm Qcoh}_G(X,\chi,W))$$
$${\rm DLFr}_G(X,\chi,W):=H^0({\rm LFr}_G(X,\chi,W))/{\rm Acycl}({\rm LFr}_G(X,\chi,W))$$
$${\rm Dlfr}_G(X,\chi,W):=H^0({\rm lfr}_G(X,\chi,W))/{\rm Acycl}({\rm lfr}_G(X,\chi,W)).$$
Objects in ${\rm Acycl}({\rm Qcoh}_G(X,\chi,W))$ are called $\textbf{acyclic}$.

\vspace{2mm}
Denote by ${\rm Acycl}^{\rm co}({\rm Qcoh}_G(X,\chi,W))$  (resp. ${\rm Acycl}^{\rm co}({\rm LFr}_G(X,\chi,W))$) the smallest thick subcategory of the triangulated category $H^0({\rm Qcoh}_G(X,\chi,W))$ (resp.  $H^0({\rm LFr}_G(X,\chi,W))$) which is closed under taking small direct sums and contain all totalizations of short exact sequences in $Z^0({\rm Qcoh}_G(X,\chi,W))$ (resp.  $Z^0({\rm LFr}_G(X,\chi,W))$).
Denote the Verdier quotients by
$${\rm D^{co}Qcoh}_G(X,\chi,W):=H^0({\rm Qcoh}_G(X,\chi,W))/{\rm Acycl^{co}}({\rm Qcoh}_G(X,\chi,W))$$
$${\rm D^{co}LFr}_G(X,\chi,W):=H^0({\rm LFr}_G(X,\chi,W))/{\rm Acycl^{co}}({\rm LFr}_G(X,\chi,W)).$$
Objects in ${\rm Acycl^{co}}({\rm Qcoh}_G(X,\chi,W))$ are called $\textbf{coacyclic}$.

If $G$ is trivial, we drop $G$ and $\chi$ from the above notation, and denote each triangulated categories by ${\rm Dcoh}(X,W)$, ${\rm DQcoh^{co}}(X,W)$, etc.

\begin{rem}
If $X$ is a regular Noetherian scheme of finite Krull dimension, then ${\rm Acycl}({\rm Qcoh}_G(X,\chi,W))$ is cocomplete, i.e. admits arbitrary direct sums (cf. \cite[Corollary 2.23]{ls}). Hence, in that case, we have 
$${\rm D^{co}Qcoh}_G(X,\chi,W)={\rm DQcoh}_G(X,\chi,W).$$
\end{rem}

\vspace{2mm}
The following lemmas ensure the existence of derived  functors between derived factorization categories.

\begin{lem}\label{inj=qcoh}
Assume that the scheme $X$ is Noetherian. The natural functor $$H^0({\rm Inj}_G(X,\chi,W))\rightarrow {\rm D^{co}Qcoh}_G(X,\chi,W)$$
is an equivalence.
\begin{proof}
Since the abelian category ${\rm Qcoh}_GX$ of $G$-equivariant quasi-coherent sheaves is a locally Noetherian Grothendieck category, it has enough injective objects, and coproducts  of injective objects are injective. Hence the result follows from \cite[Cororally 2.25]{bdfik}.
\end{proof}
\end{lem}

\begin{lem}[\cite{bfk} Proposition 3.14]\label{loc=qcoh}
Assume that $X$ is a smooth variety. Then the natural functor 
$${\rm DLFr}_G(X,\chi,W)\rightarrow{\rm DQcoh}_G(X,\chi,W)$$
is an equivalence. This equivalence induces an equivalence
$${\rm Dlfr}_G(X,\chi,W)\rightarrow{\rm Dcoh}_G(X,\chi,W).$$
\end{lem}

\subsection{Derived categories and derived factorization categories}

In this section, we recall that derived factorization categories are generalizations of bounded derived categories of coherent sheaves on schemes.

Consider  trivial $\mathbb{G}_m$-action on a scheme $X$ and an exact functor between abelian categories
$$\Upsilon: {\rm Ch}({\rm Qcoh}X)\rightarrow Z^0({\rm Qcoh}_{\mathbb{G}_m}(X,\chi_1,0))$$
given by $$\Upsilon(F^{\text{\tiny{\textbullet}}},d^{\text{\tiny{\textbullet}}}_F):=\Bigl(\bigoplus_{i\in \mathbb{Z}} F^{2i-1}(\chi_{-i})\xrightarrow{\oplus d_F^{2i-1}(\chi_{-i})}\bigoplus_{i\in \mathbb{Z}} F^{2i}(\chi_{-i})\xrightarrow{\oplus d_F^{2i}(\chi_{-i})}\bigoplus_{i\in \mathbb{Z}} F^{2i-1}(\chi_{-i+1})\Bigr).$$
Then it is easy to see that the exact functor $\Upsilon$ is an equivalence, and it induces an equivalence between triangulated categories:
$$\Upsilon: {\rm K}({\rm Qcoh}X)\rightarrow H^0({\rm Qcoh}_{\mathbb{G}_m}(X,\chi_1,0))$$
Since the triangulated equivalence preserves coacyclic objects, we obtain the following:

\begin{prop}\label{coderived}
There is an equivalence 
$$\Upsilon: {\rm D^{co}}({\rm Qcoh}X)\xrightarrow{\sim} {\rm D}^{\rm co}{\rm Qcoh}_{\mathbb{G}_m}(X,\chi_1,0)),$$
which induces an equivalence between full subcategories
$$\Upsilon: {\rm D^{b}}({\rm coh}X)\xrightarrow{\sim} {\rm D}^{}{\rm coh}_{\mathbb{G}_m}(X,\chi_1,0)).$$
\end{prop}

\begin{rem} See \cite[Appendix A]{posi1} for the definition and basic properties of coderived categories.
By \cite[Theorem 5.9.1.(b)]{posi1}, if $X$ is Noetherian, the thick subcategory ${\rm D^{b}}({\rm coh}X)$ of ${\rm D^{co}}({\rm Qcoh}X)$ is the full subcategory of compact objects. However, when $X$ is a singular variety, there is an object in ${\rm D^{b}}({\rm coh}X)$ which is not compact in the usual derived category ${\rm D}({\rm Qcoh}X)$. This is a remarkable difference between the usual derived category  and the coderived category. On the other hand, if $X$ is a smooth variety, these two kinds of derived categories are equivalent (see argument in the proof of \cite[Theorem 5.5.1.(c)]{posi1}). 
 \end{rem}

\subsection{Derived functors between derived factorization categories}
 
 We quickly review derived functors between derived factorization categories.  See, for example, \cite{ls}, \cite{bfk} or \cite{hirano} for more details.

\subsubsection{Direct images and inverse images}\hfill\\

Let $X$ and $Y$ be Noetherian schemes, and let $G$ be an affine algebraic group acting on $X$ and $Y$.
Let $f:X\rightarrow Y$ be a $G$-equivariant morphism, and choose a $\chi$-semi invariant function $W:Y\rightarrow \mathbb{A}^1$.

The morphism $f$ naturally induces a dg functor
$$f_*:{\rm Qcoh}_G(X,\chi,f^*W)\rightarrow {\rm Qcoh}_G(Y,\chi,W)$$
defined by $$f_*F:=\Bigl(f_*(F_1)\xrightarrow{f_*(\varphi_1^F)}f_*(F_0)\xrightarrow{f_*(\varphi_0^F)}f_*(F_1)(\chi)\Bigr).$$
By Lemma \ref{inj=qcoh}, we can derive the dg functor $f_*$ to obtain an exact functor 
$$\textbf{R}f_*:{\rm D^{co}Qcoh}_G(X,\chi,f^*W)\rightarrow {\rm D^{co}Qcoh}_G(Y,\chi,W).$$
If $f$ is proper morphism, it preserves factorizations whose components are coherent sheaves;
$$\textbf{R}f_*:{\rm Dcoh}_G(X,\chi,f^*W)\rightarrow {\rm Dcoh}_G(Y,\chi,W).$$

The morphism $f$ also induces another dg functor 
$$f^*:{\rm Qcoh}_G(Y,\chi,W)\rightarrow {\rm Qcoh}_G(X,\chi,f^*W),$$
defined by $$f^*E:=\Bigl(f^*(E_1)\xrightarrow{f^*(\varphi_1^E)}f^*(E_0)\xrightarrow{f^*(\varphi_0^E)}f^*(E_1)(\chi)\Bigr).$$
If $Y$ is a smooth variety, by Lemma \ref{loc=qcoh}, we have the derived functor of $f^*$
$$\textbf{L}f^*:{\rm DQcoh}_G(Y,\chi,W)\rightarrow {\rm D^{co}Qcoh}_G(X,\chi,f^*W).$$
This functor maps coherent factorizations to coherent factorizations;
$$\textbf{L}f^*:{\rm Dcoh}_G(Y,\chi,W)\rightarrow {\rm Dcoh}_G(X,\chi,f^*W).$$
It is standard that the direct image $\textbf{R}f_*$ and the inverse image $\textbf{L}f^*$ are adjoint.

\begin{rem}
If $f$ is an affine morphism, we don't need to take the derived functor, and the dg functor $f_*$ naturally defines an exact functor 
$$f_*:{\rm D^{co}Qcoh}_G(X,\chi,f^*W)\rightarrow {\rm D^{co}Qcoh}_G(Y,\chi,W).$$
Similarly, if $f$ is a flat morphism, we don't have to assume that $Y$ is a smooth variety and take the derived functor, and the dg functor $f^*$ naturally defines an exact functor 
$$f^*:{\rm D^{co}Qcoh}_G(Y,\chi,W)\rightarrow {\rm D^{co}Qcoh}_G(X,\chi,f^*W).$$
\end{rem}

\subsubsection{Tensor products}\hfill\\

Let $(X,\chi,W)^G$ be a gauged LG model, and let $V:X\rightarrow\mathbb{A}^1$ be another $\chi$-semi invariant regular function. Fix an object $F\in {\rm Qcoh}_G(X,\chi,V)$. We define a dg functor 
$$(-)\otimes F:{\rm Qcoh}_G(X,\chi,W)\rightarrow {\rm Qcoh}_G(X,\chi,W+V)$$
given by 
$$E\otimes F:=\Bigl(\bigoplus_{i=0,1}(F_{i}\otimes E_{\overline{i+1}})\xrightarrow{\varphi_1^{E\otimes F}}\bigoplus_{i=0,1}(F_{i}\otimes E_{i})(\chi^i)\xrightarrow{\varphi_0^{E\otimes F}}      \bigoplus_{i=0,1}(F_{i}\otimes E_{\overline{i+1}})(\chi)\Bigr),$$
where $\overline{n}$ is $n$ modulo 2, and 
$${\varphi_1^{E\otimes F}}=\begin{pmatrix}
\varphi^E_1\otimes1&1\otimes\varphi^F_1\\ -1\otimes\varphi^F_0&\varphi^E_0\otimes1
\end{pmatrix}
$$
and 
$$\varphi_0^{E\otimes F}=\begin{pmatrix}
\varphi^E_0\otimes1 & -(1\otimes\varphi^F_1)(\chi)\\
1\otimes\varphi^F_0 & (\varphi^E_1\otimes1)(\chi)
\end{pmatrix}.$$

By Lemma \ref{loc=qcoh}, if $X$ is a smooth variety, we have the derived functor 
$$(-)\otimes^{\textbf{L}} F:{\rm DQcoh}_G(X,\chi,W)\rightarrow {\rm DQcoh}_G(X,\chi,W+V).$$
If $F$ is a coherent factorization, this functor preserves coherent factorizations;
$$(-)\otimes^{\textbf{L}} F:{\rm Dcoh}_G(X,\chi,W)\rightarrow {\rm Dcoh}_G(X,\chi,W+V).$$
Consider a natural exact functor between abelian categories
$$\tau:{\rm Qcoh}_G(X)\rightarrow Z^0({\rm Qcoh}_G(X,\chi,0))$$
defined by $$\tau(F):=\Bigl(0\rightarrow F\rightarrow0\Bigr).$$
Then we denote by $\Sigma$ the following composition of functors
$$\Sigma:{\rm D^b}({\rm Qcoh}_GX)\xrightarrow{\tau}{\rm D^b}(Z^0({\rm Qcoh}_G(X,\chi,0)))\xrightarrow{\rm Tot}{\rm DQcoh}_G(X,\chi,0).$$
For a complex $F^{\text{\tiny{\textbullet}}}\in{\rm D^b}({\rm Qcoh}_GX)$, we define the tensor product 
$$(-)\otimes^{\textbf{L}}F^{\text{\tiny{\textbullet}}}: {\rm DQcoh}_G(X,\chi,W)\rightarrow {\rm DQcoh}_G(X,\chi,W)$$
by the following $$(-)\otimes^{\textbf{L}}F^{\text{\tiny{\textbullet}}}:=(-)\otimes^{\textbf{L}}\Sigma(F^{\text{\tiny{\textbullet}}}).$$

\begin{rem} If the components of $F$ are flat sheaves, we don't have to assume that $X$ is a smooth variety and take the derived functor, and the dg functor $(-)\otimes F$ induces an exact functor
$$(-)\otimes F:{\rm D^{co}Qcoh}_G(X,\chi,W)\rightarrow {\rm D^{co}Qcoh}_G(X,\chi,W+V).$$
Furthermore, if $F$ is a coherent factorization, the tensor product preserves coherent factorizations;
$$(-)\otimes F:{\rm Dcoh}_G(X,\chi,W)\rightarrow {\rm Dcoh}_G(X,\chi,W+V).$$
\end{rem}

\subsubsection{Integral functors}\hfill\\

We define integral functors between derived factorization categories. For simplicity, we consider  the case when $G$ is trivial. Let $X_1$ and $X_2$ be Gorenstein quasi-projective schemes, and let $W_i:X_i\rightarrow \mathbb{A}^1$ be a regular function. We denote the projection by $\pi_i:X_1\times X_2\rightarrow X_i$ for each $i=1,2$.

In order to define integral functors, we need the following lemma:

\begin{lem}[\cite{efi-posi} Corollary 2.3.e and 2.4.a]\label{fl=qcoh}
Let $(X,W)$ be a  LG model. Assume that  the scheme $X$ is a Gorenstein separated scheme of finite Krull dimension with an ample line bundle. Then the functor 
$${\rm D^{co}LFr}(X,W)\rightarrow {\rm D^{co}Qcoh}(X,W)$$
induced by the embedding of dg functor ${\rm LFr}(X,W)\rightarrow {\rm Qcoh}(X,W)$ is an equivalence.
\end{lem}

Now we define integral functors. Let $P\in{\rm D^{co}Qcoh}(X_1\times X_2,\pi_2^*W-\pi_1^*W)$ be an object.
Since $X_1$ and $X_2$ are Gorenstein, so is $X_1\times X_2$ (cf. \cite{ty}). By the above lemma, we have the derived tensor product
$$(-)\otimes^{\textbf{L}}P:{\rm D^{co}Qcoh}(X_1\times X_2,\pi_1^*W)\rightarrow{\rm D^{co}Qcoh}(X_1\times X_2,\pi_2^*W).$$

\begin{dfn}
For an object $P\in{\rm D^{co}Qcoh}(X_1\times X_2,\pi_2^*W-\pi_1^*W)$, we define the \textbf{integral functor} with respect to $P$, denoted by $\Phi_P$, as the following functor
$$\textbf{R}\pi_{2*}(\pi_1^*(-)\otimes^{\textbf{L}}P):{\rm D^{co}Qcoh}(X_1,W_1)\rightarrow{\rm D^{co}Qcoh}(X_2,W_2).$$
We call the object $P$ the \textbf{kernel} of the integral functor $\Phi_P$.
\end{dfn}

In the following, we recall that the composition of integral functors is also an integral functor. Let $X_3$ be another Gorenstein quasi-projective scheme, and let $W_3:X_3\rightarrow \mathbb{A}^1$ be a regular function. Consider the following diagram
\[\xymatrix{
&&X_1\times X_2\ar[lld]_{\pi_1}\ar[rrd]^{\pi_2}&&\\
X_1&&&&X_2\\
&&X_1\times X_2\times X_3\ar[uu]^{\pi_{12}}\ar[lld]_{\pi_{13}}\ar[rrd]^{\pi_{23}}&&\\
X_1\times X_3\ar[uu]^{q_1}\ar[rrd]_{q_3}&&&&X_2\times X_3\ar[uu]_{p_2}\ar[lld]^{p_3}\\
&&X_3&&
  }\]
where all morphisms are projections.
 For two objects
$$P\in{\rm D^{co}Qcoh}(X_1\times X_2,\pi_2^*W_2-\pi_1^*W_1)$$
$$Q\in{\rm D^{co}Qcoh}(X_2\times X_3,p_3^*W_3-p_2^*W_2),$$
we define another object
$$P\star Q:=\textbf{R}\pi_{13*}(\pi_{12}^*P\otimes^{\bf L}\pi_{23}^*Q)\in{\rm D^{co}Qcoh}(X_1\times X_3, q_3^*W_3-q_1^*W_1).$$

\begin{prop}\label{composition of integral}Notation is same as above.
The composition of integral functors 
$${\rm D^{co}Qcoh}(X_1,W_1)\xrightarrow{\Phi_P}{\rm D^{co}Qcoh}(X_2,W_2)\xrightarrow{\Phi_Q}{\rm D^{co}Qcoh}(X_3,W_3)$$
is isomorphic to the following integral functor
$${\rm D^{co}Qcoh}(X_1,W_1)\xrightarrow{\Phi_{P\star Q}}{\rm D^{co}Qcoh}(X_3,W_3).$$
\begin{proof}
The result follows from projection formula and base change formula for derived factorization categories.
\end{proof}
\end{prop}

\subsection{Support properties of factorizations}

Following \cite[Section 1.10]{efi-posi}, we consider set-theoretic supports of factorizations. In this section, $X$ is a Noetherian scheme.

\begin{dfn} Let $(X,\chi,W)^G$ be a gauged LG model, and let $Z\subset X$ be a $G$-invariant closed subset of $X$. We say that a factorization $F\in {\rm Qcoh}_G(X,\chi,W)$ is \textbf{set-theoretically supported on }$Z$ if the supports ${\rm Supp}(F_i)$ of components of $F$ are contained in $Z$.  
\end{dfn}

Denote by 
$${\rm Qcoh}_G(X,\chi,W)_Z$$
 the dg subcategory of ${\rm Qcoh}_G(X,\chi,W)$ consisting of factorizations set-theoretically supported  on $Z$.  $H^0({\rm Qcoh}_G(X,\chi,W)_Z)$ is a full triangulated subcategory of $H^0({\rm Qcoh}_G(X,\chi,W))$. Denote by ${\rm Acycl}^{\rm co}({\rm Qcoh}_G(X,\chi,W)_Z)$  the smallest thick subcategory of $H^0({\rm Qcoh}_G(X,\chi,W)_Z)$ which is closed under small direct sums and contains all totalizations of short exact sequences in $Z^0({\rm Qcoh}_G(X,\chi,W)_Z)$. Set 
$${\rm D^{co}Qcoh}_G(X,\chi,W)_Z:=H^0({\rm Qcoh}_G(X,\chi,W)_Z)/{\rm Acycl}^{\rm co}({\rm Qcoh}_G(X,\chi,W)_Z).$$
 Similarly, denote by 
 $${\rm coh}_G(X,\chi,W)_Z$$
 $${\rm Inj}_G(X,\chi,W)_Z$$
 the dg subcategory of ${\rm coh}_G(X,\chi,W)$ and ${\rm Inj}_G(X,\chi,W)$, respectively, consisting of factorizations set-theoretically supported  on $Z$. Let ${\rm Acycl}({\rm coh}_G(X,\chi,W)_Z)$  be the smallest thick subcategory of $H^0({\rm coh}_G(X,\chi,W)_Z)$ containing all totalizations of short exact sequences in $Z^0({\rm coh}_G(X,\chi,W)_Z)$, and consider the  Verdier quotient 
 $${\rm Dcoh}_G(X,\chi,W)_Z:=H^0({\rm coh}_G(X,\chi,W)_Z)/{\rm Acycl}({\rm coh}_G(X,\chi,W)_Z).$$
 
 \begin{lem}
 Let $A\in {\rm Acycl^{co}}({\rm Qcoh}_G(X,\chi,W)_Z)$ and $I\in H^0({\rm Inj}_G(X,\chi,W)_Z)$.
 Then we have 
 $${\rm Hom}_{H^0({\rm Qcoh}_G(X,\chi,W)_Z)}(A,I)=0.$$
 \begin{proof}
 Since arbitrary direct sums of  short exact sequences are exact and the totalization functor commutes with arbitrary direct sums, it is enough to show that for a short exact sequence $A^{\text{\tiny{\textbullet}}}: 0\rightarrow A^1\rightarrow A^2\rightarrow A^3\rightarrow 0$ in $Z^0({\rm Qcoh}_G(X,\chi,W)_Z)$,  we have ${\rm Hom}_{H^0({\rm Qcoh}_G(X,\chi,W)_Z)}({\rm Tot}(A^{\text{\tiny{\textbullet}}}),I)=0$. 
  This follows from a similar argument as in the proof of \cite[Lemma 2.13]{ls}.
  \end{proof}
 \end{lem}
 
 By the above Lemma, we see that every morphism  from ${\rm Acycl^{co}}({\rm Qcoh}_G(X,\chi,W)_Z)$ to ${\rm Inj}_G(X,\chi,W)_Z$ factors through the zero object in $H^0({\rm Qcoh}_G(X,\chi,W)_Z)$. Hence, by \cite[Proposition B.2]{ls}, we have the following lemma:
 
 \begin{lem}\label{right inj}
 Let $F\in H^0({\rm Qcoh}_G(X,\chi,W)_Z)$ and $I\in H^0({\rm Inj}_G(X,\chi,W)_Z)$. Then the natural map
 $${\rm Hom}_{H^0({\rm Qcoh}_G(X,\chi,W)_Z)}(F,I)\xrightarrow{\sim}{\rm Hom}_{{\rm D^{co}Qcoh}_G(X,\chi,W)_Z}(F,I)$$
 is an isomorphism.
 \end{lem}
 
 Furthermore, we have the following:

 \begin{lem}\label{inj=qcoh_z}
 The natural functor $$H^0({\rm Inj}_G(X,\chi,W)_Z)\rightarrow {\rm D^{co}Qcoh}_G(X,\chi,W)_Z$$ 
is an equivalence. 
\begin{proof}
This follows from \cite[Cororally 2.25]{bdfik}. 
\end{proof}
 \end{lem}

The following two propositions are $G$-equivariant versions of results in \cite[Section 1.10]{efi-posi}.

\begin{prop}[cf. \cite{efi-posi} Proposition 1.10]\label{support inclusion}\hfill\\
$(1)$ The natural functor $${\rm Dcoh}_G(X,\chi,W)_Z\rightarrow {\rm D^{co}Qcoh}_G(X,\chi,W)_Z$$ is fully faithful, and its image is a set of compact generators. \\
$(2)$ The natural functor $$\iota_Z:{\rm D^{co}Qcoh}_G(X,\chi,W)_Z\rightarrow {\rm D^{co}Qcoh}_G(X,\chi,W)$$ induced by the embedding of dg categories $\iota_Z:{\rm Qcoh}_G(X,\chi,W)_Z\rightarrow{\rm Qcoh}_G(X,\chi,W)$ is fully faithful.\\
$(3)$ The functor $$\iota_Z:{\rm Dcoh}_G(X,\chi,W)_Z\rightarrow {\rm Dcoh}_G(X,\chi,W)$$ induced by the embedding of dg categories $\iota_Z:{\rm coh}_G(X,\chi,W)_Z\rightarrow{\rm coh}_G(X,\chi,W)$ is fully faithful.
\begin{proof}
(1) It is enough to prove that  any morphism $F\rightarrow A$ in $ H^0({\rm Qcoh}_G(X,\chi,W)_Z)$ from  $F\in H^0({\rm coh}_G(X,\chi,W)_Z)$ to  $A\in {\rm Acycl^{co}}({\rm Qcoh}_G(X,\chi,W)_Z)$ factors through some object in ${\rm Acycl}({\rm coh}_G(X,\chi,W)_Z)$.  This follows from a similar argument  as in the proof of \cite[Lemma 2.15]{ls}. 

We show that ${\rm Dcoh}_G(X,\chi,W)_Z$ generates ${\rm D^{co}Qcoh}_G(X,\chi,W)_Z$ by using a similar discussion in the proof of \cite[Theorem 3.11.2]{posi}. By  Lemma \ref{right inj} and \ref{inj=qcoh_z}, it is enough to show that for an object $I\in H^0({\rm Inj}_G(X,\chi,W)_Z)$ if 
$${\rm Hom}_{H^0({\rm Qcoh}_G(X,\chi,W)_Z)}(F,I)=0$$
for any $F\in {\rm coh}_G(X,\chi,W)_Z$, then ${\rm id}_I:I\rightarrow I$ is homotopic to zero. Consider the partially ordered set of pairs $(M,h)$, where $M\subset I$ is a subfactorization of $I$ and $h:M\rightarrow I$ is a contracting homotopy of the embedding $i:M\hookrightarrow I$, i.e. $d(h)=i$. By Zorn's lemma, the partially ordered set contains a maximal element. Hence it suffices to show that given $(M,h)$ with $M\neq I$, there exists $(M',h')$ with $M\subsetneq M'$ and $h'|_M=h$. Take a subfactorization $M'\subset I$ such that $M\subsetneq M'$ and $M'/M\in {\rm coh}_G(X,\chi,W)_Z$. Since the components of $I$ are injective sheaves, the morphism $h:M\rightarrow I$ of degree $-1$ can be extended to a morphism $h'':M'\rightarrow I$. Denote by $i:M\hookrightarrow I$ and $i':M'\hookrightarrow I$ the embeddings. Since the map $i'-d(h'')$ is a closed degree zero morphism and vanishes on $M$, it induces a closed degree zero morphism $g:M'/M\rightarrow I$. By the assumption, $g$ is homotopic to zero, i.e. there exists a homotopy $c:M'/M\rightarrow I$ such that $d(c)=g$. Then $h'=h''+c\circ p:M'\rightarrow I$ is a contracting homotopy for $i'$ extending $h$, where $p:M'\rightarrow M'/M$ is the natural projection.

The compactness of  objects in ${\rm Dcoh}_G(X,\chi,W)_Z$ follows from Lemma \ref{right inj} and \ref{inj=qcoh_z}. 
(2) and (3) follows from Lemma \ref{inj=qcoh_z} and (1). 
\end{proof}
\end{prop}

\begin{prop}[cf. \cite{efi-posi} Theorem 1.10]\label{verdier}
Let $U:=X\setminus Z$ be the complement of $Z\subset X$, and let $j:U\rightarrow X$ be the open immersion.\\
$(1)$The restriction 
$$j^*:{\rm D^{co}Qcoh}_G(X,\chi,W)\rightarrow {\rm D^{co}Qcoh}_G(U,\chi,W|_U)$$
is the Verdier localization by the thick subcategory ${\rm D^{co}Qcoh}_G(X,\chi,W)_Z$.\\
$(2)$The restriction
$$j^*:{\rm Dcoh}_G(X,\chi,W)\rightarrow {\rm Dcoh}_G(U,\chi,W|_U)$$
is the Verdier localization by the triangulated subcategory ${\rm Dcoh}_G(X,\chi,W)_Z$. The kernel of $j^*$ is the thick envelope of ${\rm Dcoh}_G(X,\chi,W)_Z$ in ${\rm Dcoh}_G(X,\chi,W)$.
\begin{proof}
We can prove this by a standard discussion as in the proof of \cite[Theorem 1.10]{efi-posi}.

(1) Since $j^*$ has a right adjoint  $\textbf{R}j_*:{\rm D^{co}Qcoh}_G(U,\chi,W|_U)\rightarrow{\rm D^{co}Qcoh}_G(X,\chi,W)$ which is fully faithful, we see that $j^*$ is the Verdier (Bousfield) localization by its kernel which is generated by cones of adjunctions $F\rightarrow \textbf{R}j_*j^*F$ for any $F\in{\rm D^{co}Qcoh}_G(X,\chi,W)$. 

We show that ${\rm Ker}(j^*)={\rm D^{co}Qcoh}_G(X,\chi,W)_Z$. Since the inclusion ${\rm D^{co}Qcoh}_G(X,\chi,W)_Z$\\$\subset{\rm Ker}(j^*)$ is trivial, it is enough to show that the cone of the adjunction $F\rightarrow \textbf{R}j_*j^*F$, for any $F\in{\rm D^{co}Qcoh}_G(X,\chi,W)$, is contained in ${\rm D^{co}Qcoh}_G(X,\chi,W)_Z$.  By Lemma \ref{inj=qcoh}, we may take $F$ as an factorization whose components are injective quasi-coherent sheaves. Then  the adjunction comes from a closed morphism $F\rightarrow j_*j^*F$ in $Z^0({\rm Qcoh}_G(X,\chi,W))$. Since its kernel and cokernel are objects in ${\rm Qcoh}_G(X,\chi,W)_Z$, so is the cone of the adjunction by an equivariant  version of \cite[Lemma 2.7.c]{ls}.

(2) By Proposition \ref{support inclusion} (1) and \cite{nee}, we have a fully faithful functor 
$$\overline{{\rm Dcoh}_G(X,\chi,W)}/\overline{{\rm Dcoh}_G(X,\chi,W)_Z}\longrightarrow \overline{{\rm Dcoh}_G(U,\chi,W|_U)},$$
where $\overline{(-)}$ denotes the idempotent completion of $(-)$. Since every morphism $D\rightarrow E$ from $D\in{\rm Dcoh}_G(X,\chi,W)$ to $E\in \overline{{\rm Dcoh}_G(X,\chi,W)_Z}$ factors through an object in ${\rm Dcoh}_G(X,\chi,W)_Z$, the natural functor $${\rm Dcoh}_G(X,\chi,W)/{\rm Dcoh}_G(X,\chi,W)_Z\rightarrow\overline{{\rm Dcoh}_G(X,\chi,W)}/\overline{{\rm Dcoh}_G(X,\chi,W)_Z}$$ is fully faithful. Hence we see that the natural functor $${\rm Dcoh}_G(X,\chi,W)/{\rm Dcoh}_G(X,\chi,W)_Z\rightarrow{\rm Dcoh}_G(U,\chi,W|_U)$$ is also fully faithful. This functor is essentially surjective since for every $G$-equivariant coherent $F\in{\rm coh}_GU$ there exists a $G$-equivariant coherent sheaf  $\overline{F}\in{\rm coh}_GX$ such that $j^*\overline{F}\cong F$ and the coherent sheaves generate ${\rm Dcoh}_G(U,\chi,W|_U)$ by \cite[Corollary 2.29]{bdfik}
\end{proof}
\end{prop}

\subsection{Koszul factorizations}
Let $(X,\chi,W)^G$ be a gauged LG model such that $X$ is a smooth variety. Let $\mathcal{E}$ be a $G$-equivariant locally free sheaf on $X$ of rank $r$, and let 
$$s:\mathcal{E}\rightarrow \mathcal{O}_X\hspace{3mm}{\rm and}\hspace{3mm}t:\mathcal{O}_X\rightarrow\mathcal{E}(\chi)$$
be morphisms in ${\rm coh}_GX$ such that $t\circ s= W\cdot {\rm id}_{\mathcal{E}}$ and $s(\chi)\circ t=W$.
Let $Z_s\subset X$ be the zero scheme of the section $s\in \Gamma(X,\mathcal{E}^{\vee})^G$. We say that $s$ is \textbf{regular} if the codimension of $Z_s$ in $X$ equals to the rank $r$.

\begin{dfn}
We define an object $K(s,t)\in {\rm lfr}_G(X,\chi,W)$ as 
$$K(s,t):=\Bigl(K_1\xrightarrow{k_1} K_0\xrightarrow{k_0} K_1(\chi)\Bigr)$$
where 
$$K_1:=\bigoplus_{n=0}^{\lceil r/2 \rceil-1}(\bigwedge^{2n+1}\mathcal{E})(\chi^n),\hspace{4mm}{}\hspace{6mm}K_0:=\bigoplus_{n=0}^{\lfloor r/2 \rfloor}(\bigwedge^{2n}\mathcal{E})(\chi^n)$$
and $$k_i:=t\wedge(-)\oplus s\vee(-).$$
\end{dfn}

\vspace{3mm}
The following property will be necessary in section 4.1.

\begin{lem}[\cite{bfk} Lemma 3.21 and Proposition 3.20]\label{koszul lemma}\hfill\\
$(1)$ We have a natural isomorphism
$$K(s,t)^{\vee}\cong K(t^{\vee},s^{\vee}).$$
$(2)$ If $s$ is regular, we have a natural isomorphisms in ${\rm Dcoh}_G(X,\chi,W)$
$$\mathcal{O}_{Z_s}\cong K(s,t)\hspace{5mm}{\rm and}\hspace{5mm} \mathcal{O}_{Z_s}\otimes\bigwedge^{r}\mathcal{E}^{\vee}(\chi^{-1})[-r]\cong K(s,t)^{\vee},$$
where $\mathcal{O}_{Z_s}:=\Bigl(0\rightarrow\mathcal{O}_{Z_s}\rightarrow 0\Bigr)$ and $\bigwedge^{r}\mathcal{E}^{\vee}(\chi^{-1})[-r]$ is a complex in ${\rm coh}_GX$.
\end{lem}

\subsection{Restriction and Induction functors}
We construct restriction and  induction functors.
Let $G$ be an affine algebraic group acting on a scheme $X$.
Let $G\times^l X$ and $G\times ^d X$ be schemes $G\times X$ with different $G$-actions which are defined as follows;

$$
\begin{array}{ccc}
G\times G\times^l X & \stackrel{}{\longrightarrow} & G\times^l X \\
\rotatebox{90}{$\in$} & & \rotatebox{90}{$\in$} \\
(g,g',x) & \longmapsto & (gg',x)
\end{array}
$$
and
$$
\begin{array}{ccc}
G\times G\times^d X & \stackrel{}{\longrightarrow} & G\times^d X \\
\rotatebox{90}{$\in$} & & \rotatebox{90}{$\in$} \\
(g,g',x) & \longmapsto & (gg',gx).
\end{array}
$$
\vspace{2mm}
Let $\iota:X\rightarrow G\times X$ be a morphism defined by  $$X\ni x\longmapsto (1_G,x)\in G\times X.$$ We define an exact functor $\iota^*:{\rm Qcoh}_G(G\times^lX)\rightarrow{\rm Qcoh}X$ as 
$$\begin{array}{ccc}
{\rm Qcoh}_G(G\times^lX) & \stackrel{}{\longrightarrow} &{\rm Qcoh}X \\
\rotatebox{90}{$\in$} & & \rotatebox{90}{$\in$} \\
(\mathcal{F},\theta) & \longmapsto & \iota^*\mathcal{F}.
\end{array}
$$
Since $G\times^lX$ is a principal $G$-bundle over $X$, the above functor $\iota^*$ is an equivalence.

The following morphisms 
$$
\begin{array}{ccc}
\varphi: G\times^l X & \stackrel{}{\longrightarrow} & G\times^d X \\
\rotatebox{90}{$\in$} & & \rotatebox{90}{$\in$} \\
(g,x) & \longmapsto & (g,gx)
\end{array}
$$ and 
$$
\begin{array}{ccc}
\pi: G\times^d X & \stackrel{}{\longrightarrow} & X \\
\rotatebox{90}{$\in$} & & \rotatebox{90}{$\in$} \\
(g,x) & \longmapsto & x
\end{array}
$$
  are $G$-equivariant. The $G$-action $\sigma:G\times X\rightarrow X$ on $X$ is  the composition $\pi\circ\varphi$.
Since the morphism $\varphi$ is an isomorphism, the following  functors are equivalences:
$$\varphi^*:{\rm Qcoh}_G(G\times^dX)\xrightarrow{\sim}{\rm Qcoh}_G(G\times^lX)$$
$$\varphi_*:{\rm Qcoh}_G(G\times^lX)\xrightarrow{\sim}{\rm Qcoh}_G(G\times^dX).$$
Since the morphism $\pi$ is flat and affine, the following functors are exact:
 $$\pi^*:{\rm Qcoh}_G(X)\rightarrow{\rm Qcoh}_G(G\times^dX)$$
 $$\pi_*:{\rm Qcoh}_G(G\times^dX)\rightarrow{\rm Qcoh}_G(X).$$

\vspace{2mm}
\begin{dfn}\label{res ind}
We define the \textbf{restriction functor} ${\rm Res}_G:{\rm Qcoh}_GX\rightarrow{\rm Qcoh}X$ and the \textbf{induction functor} ${\rm Ind}_G:{\rm Qcoh}X\rightarrow{\rm Qcoh}_GX$ as 
$${\rm Res}_G:=\iota^*\circ\sigma^*\hspace{5mm}{\rm and}\hspace{5mm}{\rm Ind}_G:=\sigma_*\circ(\iota^*)^{-1}.$$
\end{dfn}

\vspace{2mm}
\begin{rem}
$(1)$ The restriction functor ${\rm Res}_G:{\rm Qcoh}_GX\rightarrow{\rm Qcoh}X$ is isomorphic to the forgetful functor, i.e. ${\rm Res}_G(\mathcal{F},\theta)\cong \mathcal{F}$.\\
$(2)$ Although the functor ${\rm Res}_G$ sends equivariant coherent sheaves to coherent sheaves for arbitrary group $G$, the functor ${\rm Ind}_G$ does not preserve coherent sheaves if $G$ is not a finite group.
\end{rem}

\vspace{2mm}

Since $\iota^*$ is an equivalence, the adjoint pair $\sigma^*\dashv\sigma_*$ induces the  adjoint pair
$${\rm Res}_G\dashv{\rm Ind}_{G}.$$


For a $\chi$-semi invariant regular function $W:X\rightarrow \mathbb{A}^1$, the restriction functor ${\rm Res}_G:{\rm Qcoh}_GX\rightarrow{\rm Qcoh}X$ and the induction functor ${\rm Ind}_G:{\rm Qcoh}X\rightarrow{\rm Qcoh}_GX$ induce dg functors between  factorization categories, which we denote  by the same notation; 
$${\rm Res}_G:{\rm Qcoh}_G(X,\chi,W)\rightarrow{\rm Qcoh}(X,W)$$
$${\rm Ind}_G:{\rm Qcoh}(X,W)\rightarrow{\rm Qcoh}_G(X,\chi,W).$$
These functors are also adjoint to each other;
 $${\rm Res}_G\dashv{\rm Ind}_{G}.$$
 
Since the restriction and the induction functors are exact, we obtain the following induced functors between bounded derived categories
$${\rm Res}_G:{\rm D^b}({\rm Qcoh}_GX)\rightarrow{\rm D^b}({\rm Qcoh}X)$$
$${\rm Ind}_G:{\rm D^b}({\rm Qcoh}X)\rightarrow{\rm D^b}({\rm Qcoh}_GX)$$
and the functors between coderived factorization categories
$${\rm Res}_G:{\rm D^{co}Qcoh}_G(X,\chi,W)\rightarrow{\rm D^{co}Qcoh}(X,W)$$
$${\rm Ind}_G:{\rm D^{co}Qcoh}(X,W)\rightarrow{\rm D^{co}Qcoh}_G(X,\chi,W).$$
The adjunction of the restriction and the induction functors between ${\rm Qcoh}_GX$ and ${\rm Qcoh}X$ induces the adjunction of   the above induced functors between triangulated categories;
 $${\rm Res}_G\dashv{\rm Ind}_{G}.$$


\vspace{2mm} We recall the definition of (linearly) reductiveness of algebraic groups.

\begin{dfn}
Let $H$ be an affine algebraic group over a field $K$.
\begin{itemize}
\item[$(1)$]
 $H$ is called \textbf{reductive} if the radical of $H$ is a torus.
 \item[$(2)$] $H$ is called \textbf{linearly reductive} if every rational representations of $H$ over $K$ is completely reducible.
 \end{itemize}
\end{dfn}

The following says that the above two properties are equivalent if the characteristic of the base filed is zero. 

\begin{prop}[\cite{git} Appendix A]
Let $H$ be an affine algebraic group over a field $K$ of characteristic zero.
Then $H$ is reductive if and only if it is linearly reductive.
\end{prop}


\begin{lem}\label{reductive} Assume that $G$ is linearly reductive.\\
$(1)$ The adjunction morphism $$\varepsilon:{\rm id}_{{\rm Qcoh}_GX}\rightarrow {\rm Ind}_G\circ{\rm Res}_G$$
 is a split mono i.e., there exists a  functor morphism $\eta:{\rm Ind}_G\circ{\rm Res}_G\rightarrow {\rm id}_{{\rm Qcoh}_GX}$ such that $\eta\circ\varepsilon={\rm id}$. The adjunction morphism $${\rm id}_{{\rm Qcoh}_G(X,\chi,W)}\rightarrow {\rm Ind}_G\circ{\rm Res}_G$$
 is also a split mono.\\
$(2)$ The restriction functors 
$${\rm Res}_G:{\rm D^b}({\rm Qcoh}_GX)\rightarrow{\rm D^b}({\rm Qcoh}X)$$ 
$${\rm Res}_G:{\rm D^{co}Qcoh}_G(X,\chi,W)\rightarrow{\rm D^{co}Qcoh}(X,W)$$ are faithful.

\begin{proof}
(1) This follows from the argument in the proof of \cite[Lemma 4.56]{hirano}.\\
(2) We will prove that the upper functor ${\rm Res}_G:{\rm D^b}({\rm Qcoh}_GX)\rightarrow{\rm D^b}({\rm Qcoh}X)$ is faithful; the proof of the faithfulness of the lower functor is similar. The functor morphism $\eta:{\rm Ind}_G\circ{\rm Res}_G\rightarrow {\rm id}_{{\rm Qcoh}_GX}$ constructed in (1) naturally induces the functor morphism $\overline{\eta}:{\rm Ind}_G\circ{\rm Res}_G\rightarrow {\rm id}_{{\rm D^b}({\rm Qcoh}_GX)}$ such that the composition with the adjunction morphism $${\rm id}_{{\rm D^b}({\rm Qcoh}_GX)}\rightarrow {\rm Ind}_G\circ{\rm Res}_G\xrightarrow{\overline{\eta}} {\rm id}_{{\rm D^b}({\rm Qcoh}_GX)}$$ is the identity. Hence any morphism $f$ in ${\rm D^b}({\rm Qcoh}_GX)$  factors through ${\rm Ind}_G\circ{\rm Res}_G(f)$, and so $f=0$ if ${\rm Res}(f)=0$.  
\end{proof}
\end{lem}


\section{Relative singularity categories}
Relative singularity categories are introduced in \cite{efi-posi}, and it is shown that derived factorization categories (with some conditions on regular functions) are equivalent to relative singularity categories. In this section, we recall the definition and  properties of relative singularity categories.

\subsection{Triangulated categories of relative singularities}
 
Let $X$ be a  quasi-projective scheme, and let $G$ be an affine algebraic group acting on $X$. Throughout this section, we assume that $X$ {\it has a $G$-equivariant ample line bundle}. If $X$ is normal, this condition is satisfied by \cite[Lemma 2.10]{thomason}. The equivariant triangulated category of singularities ${\rm D}^{\rm sg}_G(X)$ of $X$ is defined as the Verdier quotient of ${\rm D^b}({\rm coh}_GX)$ by the thick subcategory ${\it Perf}_G(X)$ of equivariant perfect complexes. Following \cite{orlov2},  we consider a larger category ${\rm D}^{\rm cosg}_G(X)$ defined as the Verdier quotient of ${\rm D^b}({\rm Qcoh}_GX)$ by the thick subcategory ${\it Lfr}_G(X)$ of complexes which is  quasi-isomorphic to  bounded complexes of equivariant locally free  sheaves (not necessarily of finite ranks). If $G$ is trivial, we denote the singularity categories by ${\rm D}^{\rm cosg}(X)$ or ${\rm D}^{\rm sg}(X)$.

We recall relative singularity categories following \cite{efi-posi}. Let  $Z\subset X$ be a $G$-invariant closed subscheme of $X$ such that $\mathcal{O}_{Z}$  {\it has finite $G$-flat dimension as an $\mathcal{O}_{X}$-module} i.e., the $G$-equivariant sheaf $\mathcal{O}_Z\in{\rm coh}_G(X)$ has a finite resolution $F^{\text{\tiny{\textbullet}}}\rightarrow\mathcal{O}_Z$ of $G$-equivariant flat sheaves on $X$. Under the assumption, we have the derived inverse image $\mathbf{L}i^*:{\rm D^b}({\rm Qcoh}_GX)\rightarrow {\rm D^b}({\rm Qcoh}_GZ)$ between bounded derived categories for the closed immersion $i:Z \rightarrow X$. This functor preserves complexes of coherent sheaves; $\mathbf{L}i^*:{\rm D^b}({\rm coh}X)\rightarrow {\rm D^b}({\rm coh}Z)$. 

\begin{dfn}[\cite{efi-posi} Section 2.1]
We define  the following Verdier quotients 
$${\rm D}_G^{\rm cosg}(Z/X):={\rm D^b}({\rm Qcoh}_GZ)/\langle {\rm Im}(\mathbf{L}i^*:{\rm D^b}({\rm Qcoh}_GX)\rightarrow {\rm D^b}({\rm Qcoh}_GZ))\rangle^{\oplus} $$
$${\rm D}_G^{\rm sg}(Z/X):={\rm D^b}({\rm coh}_GZ)/\langle {\rm Im}(\mathbf{L}i^*:{\rm D^b}({\rm coh}_GX)\rightarrow {\rm D^b}({\rm coh}_GZ))\rangle, $$
where $\langle - \rangle$ (resp. $\langle - \rangle^{\oplus}$) denotes the smallest thick subcategory containing objects in $(-)$ (resp. and closed under infinite direct sums which exist in ${\rm D^b}({\rm Qcoh}_GZ)$). The quotient category ${\rm D}^{\rm sg}_{G}(Z/X)$ is called the \textbf{equivariant triangulated category of singularities of} $Z$ \textbf{relative to} $X$. If $G$ is trivial, we denote the categories defined above by  ${\rm D}^{\rm cosg}(Z/X)$ or ${\rm D}^{\rm sg}(Z/X)$.
\end{dfn}

\begin{prop}\label{verdier sing}
Assume that $G$ is reductive. We have  natural Verdier localizations by thick subcategories
$$\pi^{\rm co}:{\rm D}^{\rm cosg}_G(Z)\rightarrow{\rm D}^{\rm cosg}_G(Z/X)$$ 
$$\pi:{\rm D}^{\rm sg}_G(Z)\rightarrow{\rm D}^{\rm sg}_G(Z/X).$$
\begin{proof}
It is enough to show that ${\it Lfr}_G(Z)\subset\langle {\rm Im}(\mathbf{L}i^*:{\rm D^b}({\rm Qcoh}_GX)\rightarrow {\rm D^b}({\rm Qcoh}_GZ))\rangle^{\oplus}$ and 
${\it Perf}_G(Z)\subset \langle {\rm Im}(\mathbf{L}i^*:{\rm D^b}({\rm coh}_GX)\rightarrow {\rm D^b}({\rm coh}_GZ))\rangle$.
These inclusions follow from the assumption that $X$ has a $G$-equivariant ample line bundle $L$.
The proofs of the inclusions are similar, and we  prove the only former inclusion.
It is enough to show that any $G$-equivariant locally free sheaf $E$ on $Z$ is a direct summand of a bounded complex whose terms are direct sums of invertible sheaves of the form $i^*L^{\otimes n}$. By \cite[Lemma 1.4]{thomason}, there is a bounded above locally free resolution $E^{\text{\tiny{\textbullet}}}\xrightarrow{\sim} E$  whose terms are as above. For any $n>0$, we have the following triangle in ${\rm D^b}({\rm Qcoh}_GZ)$
$$\sigma^{\geq-n}E^{\text{\tiny{\textbullet}}}\rightarrow E\rightarrow H^{-n}(\sigma^{\geq-n}E^{\text{\tiny{\textbullet}}})[n+1]\rightarrow \sigma^{\geq-n}E^{\text{\tiny{\textbullet}}}[1],$$
where $\sigma^{\geq-n}$ denotes the brutal truncation. If we choose a sufficiently large $n\gg0$, we have 
$${\rm Hom}_{{\rm D^b}({\rm Qcoh}_GZ)}(E,H^{-n}(\sigma^{\geq-n}E^{\text{\tiny{\textbullet}}})[n+1])=0$$
by \cite[Lemma 1.12]{orlov2}, since the restriction functor ${\rm Res}_G:{\rm D^b}({\rm Qcoh}_GZ)\rightarrow{\rm D^b}({\rm Qcoh}Z)$ is faithful by  Lemma \ref{reductive} (2). Hence the above triangle splits, and $E$ is a direct summand of the complex $\sigma^{\geq-n}E^{\text{\tiny{\textbullet}}}$.
\end{proof}
\end{prop}

\begin{rem}\label{rel equal usual}
Note that, if $X$ is regular, then the thick subcategory $\langle {\rm Im}(\mathbf{L}i^*)\rangle\subset  {\rm D^b}({\rm coh}_GZ)$ coincides with its thick subcategory ${\it Perf}_G(Z)$ of equivariant perfect complexes of $Z$. Hence the quotient category  ${\rm D}_{G}^{\rm sg}(Z/X)$ is same as ${\rm D}_G^{\rm sg}(Z)$. Similarly, ${\rm D}_{G}^{\rm cosg}(Z/X)$ is also same as ${\rm D}_G^{\rm cosg}(Z)$ when $X$ is regular.
\end{rem}

The exact functors ${\rm Res}_G:{\rm Qcoh}_GZ\rightarrow {\rm Qcoh}Z$ and  ${\rm Ind}_G:{\rm Qcoh}Z\rightarrow {\rm Qcoh}_GZ$, defined in Definition \ref{res ind}, induce functors between relative singularity categories
$${\rm Res}_G:{\rm D^{cosg}_G}(Z/X)\rightarrow {\rm D^{cosg}}(Z/X)$$
$${\rm Ind}_G:{\rm D^{cosg}}(Z/X)\rightarrow {\rm D^{cosg}_G}(Z/X).$$
We need the following lemma in the proof of the main result.

\begin{lem}\label{faithful lemma}
Assume that $G$ is reductive. Then the restriction functor 
$${\rm Res}_G:{\rm D^{cosg}_G}(Z/X)\rightarrow {\rm D^{cosg}}(Z/X)$$ 
is  faithful.
\begin{proof}
This follows from a similar argument as in the proof of Lemma \ref{reductive} (2).
\end{proof}
\end{lem}

\subsection{Direct images and inverse images in relative singularity categories}
 Let $X_1$ and $X_2$ be quasi-projective schemes with actions of an affine algebraic group $G$. Assume that $X_1$ and $X_2$ have $G$-equivariant ample line bundles. Let  $\tilde{f}:X_2\rightarrow X_1$ be a $G$-equivariant morphism. Let $Z_1$ be a $G$-invariant closed subscheme of $X_1$ such that $\mathcal{O}_{Z_1}$ has finite $G$-flat dimension as a $\mathcal{O}_{X_1}$-module, and let $Z_2$ be the fiber product $Z_1\times_{X_1}X_2$. Denote by $f$ the restriction $\tilde{f}|_{Z_2}:Z_2\rightarrow Z_1$ of $\tilde{f}$ to $Z_2$. We assume that the cartesian square 
 \[\xymatrix{
Z_2\ar[rr]^{f}\ar[d]_{}&&Z_1\ar[d]^{}\\
X_2\ar[rr]^{\tilde{f}}&&X_1
}\]
is {\it exact} in the sense of \cite{kuz}. Then, $\mathcal{O}_{Z_2}$ also has finite $G$-flat dimension as a $\mathcal{O}_{X_2}$-module.
   Furthermore, we assume that  $\tilde{f}$ {\it has finite G-flat dimension}, i.e.  the derived inverse image $\textbf{L}\tilde{f}^*:{\rm D}^{-}({\rm Qcoh}_GX_1)\rightarrow {\rm D}^{-}({\rm Qcoh}_GX_2)$ maps ${\rm D^b}({\rm Qcoh}_GX_1)$ to ${\rm D^b}({\rm Qcoh}_GX_2)$. Then $f$ also has finite $G$-flat dimension. 
   
   In the above setting, the derived inverse image $\textbf{L}f^*:{\rm D^b}({\rm Qcoh}_GZ_1)\rightarrow {\rm D^b}({\rm Qcoh}_GZ_2)$ induces exact functors
 $$f^{\circ}:{\rm D}^{\rm cosg}_G(Z_1/X_1)\rightarrow {\rm D}^{\rm cosg}_G(Z_2/X_2)$$
 $$f^{\circ}:{\rm D}^{\rm sg}_G(Z_1/X_1)\rightarrow {\rm D}^{\rm sg}_G(Z_2/X_2),$$
 and the derived direct image $\textbf{R}f_*:{\rm D^b}({\rm Qcoh}_GZ_2)\rightarrow{\rm D^b}({\rm Qcoh}_GZ_1)$ induces a right adjoint functor of $f^{\circ}:{\rm D}^{\rm cosg}_G(Z_1/X_1)\rightarrow {\rm D}^{\rm cosg}_G(Z_2/X_2)$
  $$f_{\circ}:{\rm D}^{\rm cosg}_G(Z_2/X_2)\rightarrow {\rm D}^{\rm cosg}_G(Z_1/X_1).$$
If $f$ is a proper morphism, the direct image $\textbf{R}f_*:{\rm D^b}({\rm coh}_GZ_2)\rightarrow{\rm D^b}({\rm coh}_GZ_1)$ between bounded complexes of coherent sheaves induces a right adjoint functor   $$f_{\circ}:{\rm D}^{\rm sg}_G(Z_2/X_2)\rightarrow {\rm D}^{\rm sg}_G(Z_1/X_1)$$ 
  of $f^{\circ}:{\rm D}^{\rm sg}_G(Z_1/X_1)\rightarrow {\rm D}^{\rm sg}_G(Z_2/X_2)$.

Let $X$ be a quasi-projective scheme with an action of an affine algebraic group $G$, and let $U\subset X$ be a $G$-invariant open subscheme.  Let $Z\subset X$ be a $G$-invariant closed subscheme such that $\mathcal{O}_Z$ has finite $G$-flat dimension, and consider the fiber product $U_Z:=Z\times_X U$. Denote by $\tilde{l}:U\rightarrow X$ and $l:U_Z\rightarrow Z$ the open immersions. Then we have the following exact cartesian square:
\[\xymatrix{
U_Z\ar[rr]^{l}\ar[d]_{}&&Z\ar[d]^{}\\
U\ar[rr]^{\tilde{l}}&&X.
}\]

\begin{lem}\label{local verdier} The inverse image $$l^{\circ}:{\rm D}_G^{\rm cosg}(Z/X)\rightarrow {\rm D}_G^{\rm cosg}(U_Z/U)$$
 is a Verdier localization by the kernel of $l^{\circ}$.
\begin{proof}
The direct image $\textbf{R}l_*: {\rm D^b}({\rm Qcoh}U_Z)\rightarrow{\rm D^b}({\rm Qcoh}Z)$ is fully faithful and right adjoint to the inverse image $l^*:{\rm D^b}({\rm Qcoh}Z)\rightarrow{\rm D^b}({\rm Qcoh}U_Z)$. By \cite[Lemma 1.1]{orlov1}, the direct image functor $l_{\circ}:{\rm D}_G^{\rm cosg}(Z/X)\rightarrow {\rm D}_G^{\rm cosg}(U_Z/U)$ is fully faithful. Hence, $l^{\circ}$ admits a right adjoint functor which is fully faithful, and this implies the result.
\end{proof}
\end{lem}

\subsection{Relative singularity categories and derived factorization categories}
\vspace{3mm}
 In this section, $X$ and $G$ are the same as in section 3.1, and we assume that $G$ is {\it reductive}. Let $\chi:G\rightarrow \mathbb{G}_m$ be a character of $G$, and let $W: X\rightarrow \mathbb{A}^1$ be a $\chi$-semi invariant regular function. In this section, we assume that the corresponding $G$-invariant section $W:\mathcal{O}_X\rightarrow \mathcal{O}(\chi)$ is injective. For example, if $W$ is flat, this condition is satisfied. Denote by $X_0$ the fiber of $W$ over $0\in \mathbb{A}^1$, and let $i:X_0\rightarrow X$ be the closed immersion. We have an exact functor $\tau:{\rm Qcoh}_GX_0\rightarrow Z^0({\rm Qcoh}_G(X,\chi,W))$ defined by 
$$\tau(F):=\Bigl( 0\rightarrow i_*(F)\rightarrow 0\Bigr).$$
We define a natural functor 
$$\Upsilon:{\rm D^b}({\rm Qcoh}_GX_0)\rightarrow {\rm D^{co}Qcoh}_G(X,\chi,W)$$
as the composition of functors 
$${\rm D^b}({\rm Qcoh}_GX_0)\xrightarrow{\tau}{\rm D^b}(Z^0({\rm Qcoh}_G(X,\chi,W)))\xrightarrow{\rm Tot} {\rm D^{co}Qcoh}_G(X,\chi,W).$$
The functor $\Upsilon$ annihilates the thick category $\langle {\rm Im}(\textbf{L}i^*)\rangle^{\oplus}\subset{\rm D^b}({\rm Qcoh}_GX_0)$, since its nonequivariant functor $\Upsilon:{\rm D^b}({\rm Qcoh}X_0)\rightarrow {\rm D^{co}Qcoh}(X,W)$ annihilates ${\rm Res}_{G}(\langle {\rm Im}(\textbf{L}i^*)\rangle^{\oplus})$ (see the proof of \cite[Theorem 2.7 and Theorem 2.8]{efi-posi}) and the restriction functor ${\rm Res}_G:{\rm D^{co}Qcoh}_G(X,\chi,W)\rightarrow{\rm D^{co}Qcoh}(X,W)$ is faithful. Hence it induces an exact functor 
$$\Upsilon:{\rm D}^{\rm cosg}_G(X_0/X)\rightarrow {\rm D^{co}Qcoh}_G(X,\chi,W).$$
Similarly, we have the following exact functor  
$$\Upsilon:{\rm D}^{\rm sg}_G(X_0/X)\rightarrow {\rm Dcoh}_G(X,\chi,W),$$
and the following diagram is commutative;
 \[\xymatrix{
{\rm D}^{\rm cosg}_G(X_0/X)\ar[rr]^{\Upsilon}&&{\rm D^{co}Qcoh}_G(X,\chi,W)\\
{\rm D}^{\rm sg}_G(X_0/X)\ar[rr]^{\Upsilon}\ar[u]&&{\rm Dcoh}_G(X,\chi,W),\ar[u]
}\]
where the vertical arrows are natural inclusion functors (which are fully faithful).


\begin{thm}[cf. \cite{efi-posi} Theorem 2.7, Theorem 2.8.]\label{relative sing}
The functors $$\Upsilon:{\rm D}^{\rm cosg}_G(X_0/X)\rightarrow {\rm D^{co}Qcoh}_G(X,\chi,W)$$ 
$$\Upsilon:{\rm D}^{\rm sg}_G(X_0/X)\rightarrow {\rm Dcoh}_G(X,\chi,W)$$ 
are equivalences.
\end{thm}

In order to prove the above theorem, we need to construct the quasi-inverse of $\Upsilon$.
We say that a $G$-equivariant quasi-coherent sheaf $F\in{\rm Qcoh}_GX$ is \textbf{$W$-flat}, if the morphism of sheaves $W:F\rightarrow F\otimes L$ is injective. Denote by ${\rm Flat}^W_G(X,\chi,W)$  the dg full subcategory of ${\rm Qcoh}_G(X,\chi,W)$  consisting of factorizations whose components are $W$-flat. Then $H^0({\rm Flat}^W_G(X,\chi,W))$   is a full triangulated subcategory of $H^0({\rm Qcoh}_G(X,\chi,W))$.
Denote by ${\rm Acycl^{co}}({\rm Flat}_G^W(X,\chi,W))$  the smallest thick subcategory of $H^0({\rm Flat}^W_G(X,\chi,W))$  containing all totalizations of short exact sequences in the exact category $Z^0({\rm Flat}^W_G(X,\chi,W))$. Consider  the corresponding Verdier quotients
$${\rm D^{co}Flat}_G^W(X,\chi,W):=H^0({\rm Flat}^W_G(X,\chi,W))/{\rm Acycl^{co}}({\rm Flat}_G^W(X,\chi,W)).$$
The restriction functor ${\rm Res}_G:{\rm Qcoh}_G(X,\chi,W)\rightarrow{\rm Qcoh}(X,W)$ and the induction functor ${\rm Ind}_G:{\rm Qcoh}(X,W)\rightarrow{\rm Qcoh}_G(X,\chi,W)$  preserve factorizations whose components are $W$-flat sheaves since ${\rm Res}_G:{\rm Qcoh}_G X\rightarrow {\rm Qcoh}X$ and ${\rm Ind}_G:{\rm Qcoh} X\rightarrow {\rm Qcoh}_GX$ are exact functors. Hence the restriction and the induction functors induce the following functors
$${\rm Res}_G:{\rm D^{co}Flat}_G^W(X,\chi,W)\rightarrow{\rm D^{co}Flat}^W(X,W)$$ 
$${\rm Ind}_G:{\rm D^{co}Flat}^W(X,W)\rightarrow{\rm D^{co}Flat}_G^W(X,\chi,W),$$ 
and these functors are adjoint to each other;
  $${\rm Res}_G\dashv{\rm Ind}_{G}.$$

\begin{lem}\label{w-flat}
The natural functor
$${\rm D^{co}Flat}_G^W(X,\chi,W)\rightarrow{\rm D^{co}Qcoh}_G(X,\chi,W)$$
is an equivalence.
\begin{proof}
At first, we prove that the functor is essentially surjective. Let $F\in {\rm D^{co}Qcoh}_G(X,\chi,W)$ be an object.
Since $X$ has a $G$-equivariant ample line bundle, there are $G$-equivariant locally free sheaf $E_i$  and a surjective morphism $p_i:E_i\rightarrow F_i$ in ${\rm Qcoh}_G X$ for each $i=0,1$. Let  $E\in{\rm Qcoh}_G(X,\chi,W)$ be the factorization of the following form
$$E:=\Bigl(E_1\oplus E_0\xrightarrow{W\oplus {\rm id}_{E_0}} E_1(\chi)\oplus E_0\xrightarrow{{\rm id}_{E_1(\chi)}\oplus W}E_1(\chi)\oplus E_0(\chi)\Bigr).$$
Then $p_1$ and $p_0$ define a natural surjective morphism $p:E\rightarrow F$ in $Z^0({\rm Qcoh}_G(X,\chi,W))$. The kernel $K:={\rm Ker}(p)$ of $p$ is in $Z^0({\rm Flat}_G^W(X,\chi,W))$ since the components of $K$ are subsheaves of $W$-flat sheaves. Hence the totalization ${\rm Tot}(C^{\text{\tiny{\textbullet}}})$ of the complex 
$$C^{\text{\tiny{\textbullet}}}:\cdots \rightarrow 0\rightarrow K\hookrightarrow E\rightarrow0\rightarrow \cdots$$ with the cohomological degree of $E$  zero is in ${\rm D^{co}Flat}_G^W(X,\chi,W)$, and we see that the natural  morphism  ${\rm Tot}(C^{\text{\tiny{\textbullet}}})\rightarrow F$ induced by $p$ is an isomorphism  in ${\rm D^{co}Qcoh}_G(X,\chi,W)$.

To show the functor ${\rm D^{co}Flat}_G^W(X,\chi,W)\rightarrow{\rm D^{co}Qcoh}_G(X,\chi,W)$ is fully faithful, it suffices to prove that for any morphism $f:E\rightarrow F$ in $H^0({\rm Qcoh}_G(X,\chi,W))$ with $F\in H^0({\rm Flat}^W_G(X,\chi,W))$ and the cone of $f$ in ${\rm Acycl^{co}}({\rm Qcoh}_G(X,\chi,W))$, there exists a morphism $g:F'\rightarrow E$ with $F'\in H^0({\rm Flat}^W_G(X,\chi,W))$ such that the cone of $f\circ g$ is in ${\rm Acycl^{co}}({\rm Flat}_G^W(X,\chi,W))$ (see \cite[Proposition B.2. $({\it ff1})^{\rm op}$]{ls}). By the above argument in the previous paragraph, we can find a morphism  $g:F'\rightarrow E$ with $F'\in H^0({\rm Flat}^W_G(X,\chi,W))$ such that the cone  of $g$ is in ${\rm Acycl^{co}}({\rm Qcoh}_G(X,\chi,W))$, and then the cone of $f\circ g$ is in $H^0({\rm Flat}^W_G(X,\chi,W))\cap{\rm Acycl^{co}}({\rm Qcoh}_G(X,\chi,W))$. Hence it is enough to show that $$H^0({\rm Flat}^W_G(X,\chi,W))\cap{\rm Acycl^{co}}({\rm Qcoh}_G(X,\chi,W))\subseteq{\rm Acycl^{co}}({\rm Flat}_G^W(X,\chi,W)).$$
For this, let $A\in H^0({\rm Flat}^W_G(X,\chi,W))\cap{\rm Acycl^{co}}({\rm Qcoh}_G(X,\chi,W))$ be an object. We already know that ${\rm Res}_G(A)\in {\rm Acycl^{co}}({\rm Flat}^W(X,W))$ by \cite[Corollary 2.6 (a)]{efi-posi}. Note that the restriction functor ${\rm Res}_G:{\rm D^{co}Flat}_G^W(X,\chi,W)\rightarrow{\rm D^{co}Flat}^W(X,W)$ is faithful by a similar argument as in the proof of Lemma \ref{reductive} (2). Hence the fact that  ${\rm Res}_G(A)\in {\rm Acycl^{co}}({\rm Flat}^W(X,W))$ implies that  $A\in {\rm Acycl^{co}}({\rm Flat}_G^W(X,\chi,W))$.
\end{proof}
\end{lem}

\vspace{2mm}
For an object  $F=\Bigl(F_1\xrightarrow{\varphi_1^F}F_0\xrightarrow{\varphi_0^F}F_1(\chi)\Bigr)\in Z^0({\rm Flat}^W_G(X,\chi,W))$, define an object $\Xi(F)\in{\rm D}_G^{\rm cosg}(X_0/X)$ by $\Xi(F):={\rm Cok}(\varphi_1^F)$. It is easy to see that this defines the following exact functor
$$\Xi:H^0({\rm Flat}^W_G(X,\chi,W))\rightarrow{\rm D}_G^{\rm cosg}(X_0/X).$$ 
If $G $ is trivial, this exact functor annihilates  ${\rm Acycl^{co}}({\rm Flat}^W(X,W))$ by \cite[Theorem 2.7, 2.8]{efi-posi}. Hence, since  ${\rm Res}_G:{\rm D}_G^{\rm cosg}(X_0/X)\rightarrow{\rm D}^{\rm cosg}(X_0/X)$ is faithful, we obtain  the exact functor 
 $\Xi:{\rm D^{co}Flat}_G^W(X,\chi,W)\rightarrow {\rm D}_G^{\rm cosg}(X_0/X).$ By Lemma \ref{w-flat}, we have 
the left derived functor of $\Xi$;
$$\textbf{L}\Xi: {\rm D^{co}Qcoh}_G(X,\chi,W)\rightarrow {\rm D}^{\rm cosg}_G(X_0/X).$$

\vspace{5mm}

{\it Proof of Theorem \ref{relative sing}:}
 We will show that the functors $\Upsilon$ and $\textbf{L}\Xi$ are mutually inverse. 
 Let $E\in {\rm D^{co}Qcoh}_G(X,\chi,W)$ be an object. By Lemma \ref{w-flat} we may assume that $E\in{\rm D^{co}Flat}^W_G(X,\chi,W)$.
Then $$\Upsilon\textbf{L}\Xi(E)\cong \Upsilon\Xi(E)=\Bigl(0\rightarrow{\rm Cok}(\varphi_1^E)\rightarrow0\Bigr),$$ and the  surjective morphism $E_0\twoheadrightarrow{\rm Cok}(\varphi_1^E)$ induces the  natural surjective morphism $\phi_E:E\rightarrow \Upsilon\Xi(E)$ in $Z^0({\rm Qcoh}_G(X,\chi,W))$. Since the kernel of $\phi_E$ is the factorization $\Bigl(E_1=E_1\xrightarrow{W}E_1(\chi)\Bigr)$ and it is isomorphic to the zero object in $H^0({\rm Flat}_G^W(X,\chi,W))$, the morphism $\phi_E:E\rightarrow \Upsilon\Xi(E)$ is an isomorphism in ${\rm D^{co}Qcoh}_G(X,\chi,W)$. 
It is easy to see that the isomorphisms  $\phi_{(-)}$ define an isomorphism of functors $$\phi:{\rm id}_{{\rm D^{co}Qcoh}_G(X,\chi,W)}\xrightarrow{\sim}\Upsilon\textbf{L}\Xi.$$

Let $F\in{\rm D}^{\rm cosg}_G(X_0/X)$ be an object. Then we may assume that $F\in{\rm Qcoh}_GX_0$. Take a surjective morphism $p:P\twoheadrightarrow i_*F$ with $P$ locally free. Set $K:={\rm Ker}(p)\in {\rm Qcoh}_GX$ and $Q:=\bigl(K\xrightarrow{i} P\xrightarrow{W}K(\chi)\bigr)\in {\rm Qcoh}_G(X\chi,W)$, where $i:K\rightarrow P$ is the natural inclusion. Consider the natural surjective morphism $\pi:Q\rightarrow\bigl(0\rightarrow i_*F\rightarrow0\bigr)$
 in $Z^0({\rm Qcoh}_G(X,\chi,W))$. Then the kernel of $\pi$ is the factorization $\bigl(K=K\xrightarrow{W}K(\chi)\bigr)$, and it is isomorphic to the zero object in $H^0({\rm Qcoh}_G(X,\chi,W))$. Hence $\pi$ is an isomorphism in ${\rm D^{co}Qcoh}_G(X,\chi,W)$, and so we have a natural isomorphism $\psi_F:\textbf{L}\Xi \Upsilon(F)\xrightarrow{\sim}F$ in ${\rm D}^{\rm cosg}_G(X_0/X)$
 defined as the composition $\textbf{L}\Xi \Upsilon(F)\xrightarrow{\sim}\Xi \Upsilon(Q)={\rm Cok}(i)=F$.
 We need to show that the isomorphisms $\psi_{(-)}$ are functorial in $(-)$. Since the restriction functor ${\rm Res}_G$ is isomorphic to the forgetful functor ${\rm Forg}_G$, we have a natural isomorphism of functors $\sigma:{\rm Res}_G\textbf{L}\Xi \Upsilon\xrightarrow{\sim}\textbf{L}\Xi \Upsilon{\rm Res}_G$ defined by the composition $${\rm Res}_G\textbf{L}\Xi \Upsilon\xrightarrow{\sim}{\rm Forg}_G\textbf{L}\Xi \Upsilon=\textbf{L}\Xi \Upsilon{\rm Forg}_G\xrightarrow{\sim}\textbf{L}\Xi \Upsilon{\rm Res}_G,$$
  and the following diagram is commutative
  \[\xymatrix{
{\rm Res}_G\textbf{L}\Xi \Upsilon(F)\ar[rr]^{{\rm Res}_G(\psi_F)}\ar[rd]_{\sigma_F}&&{\rm Res}_G(F)\\
&\textbf{L}\Xi \Upsilon{\rm Res}_G(F)\ar[ru]_{\psi_{{\rm Res_G}(F)}}&
}\]
Hence we see that the isomorphisms $\psi_{(-)}$ are functorial by the fact that   the isomorphisms $\psi_{(-)}$ are  functorial if $G$ is trivial and that the functor ${\rm Res}_G$ is faithful. This completes the proof of the former equivalence.

The latter equivalence follows from \cite[Remark 2.7]{efi-posi}, which is a generalized result of \cite[Theorem 2.7]{efi-posi}. \qed


\section{Derived Kn\"orrer Periodicity}

Let $X$ be a smooth quasi-projective variety, and let $G$ be a reductive affine algebraic group acting on $X$. Let $\mathcal{E}$ be a $G$-equivariant locally free sheaf of rank $r$, and let $s\in \Gamma(X,\mathcal{E}^{\vee})^G$ be a $G$-invariant section of $\mathcal{E}^{\vee}$. Denote by $Z\subset X$ the zero scheme of $s$. We assume that $s$ is {\it regular}, i.e. the codimension of $Z$ in $X$ is $r$. Let
$$\rm{V}(\mathcal{E}(\chi)):=\underline{\rm Spec}(Sym(\mathcal{E}(\chi)^{\vee}))$$
be a vector bundle over $X$ with the $G$-action induced by the equivariant structure of the locally free sheaf $\mathcal{E}(\chi)$. Denote by ${\rm V}(\mathcal{E}(\chi))|_Z$ the restriction of the vector bundle $\rm{V}(\mathcal{E}(\chi))$ to $Z$. Let $j:Z\hookrightarrow X$ and $i:{\rm V}(\mathcal{E}(\chi))|_Z\hookrightarrow\rm{V}(\mathcal{E}(\chi))$ be the closed immersions, and let $q:\rm{V}(\mathcal{E}(\chi))\rightarrow X$ and $p:{\rm V}(\mathcal{E}(\chi))|_Z\rightarrow Z$ be the projections. Now we have the following commutative diagram:
\[\xymatrix{
{\rm V}(\mathcal{E}(\chi))|_Z\ar[r]^{i}\ar[d]_{p}&\rm{V}(\mathcal{E}(\chi))\ar[d]^{q}\\
Z\ar[r]^{j}&X.
}\]
The invariant section $s$ induces a $\chi$-semi invariant regular function $$Q_s:\rm{V}(\mathcal{E}(\chi))\rightarrow\mathbb{A}^1.$$
Let $W:X\rightarrow\mathbb{A}^1$ be a $\chi$-semi invariant regular function on $X$. The function $W$ induces  $\chi$-semi invariant functions on $Z$, ${\rm V}(\mathcal{E}(\chi))$ and  ${\rm V}(\mathcal{E}(\chi))|_Z$, which we denote by the same notation $W$ (by abuse of notation).
Since the inverse image $p^*$ and the direct image $i_*$ are exact and commutative with arbitrary direct sums as functors between categories of quasi-coherent sheaves, these induce (underived) functors $$p^* : {\rm D^{co}Qcoh}_G(Z,\chi,W)\rightarrow{\rm D^{co}Qcoh}_G({\rm V}(\mathcal{E}(\chi))|_Z,\chi,W)$$
$$i_*:{\rm D^{co}Qcoh}_G({\rm V}(\mathcal{E}(\chi))|_Z,\chi,W)\rightarrow{\rm D^{co}Qcoh}_G({\rm V}(\mathcal{E}(\chi)),\chi,W+Q_s).$$
Restricting the composition $i_*p^*:{\rm D^{co}Qcoh}_G(Z,\chi,W)\rightarrow{\rm D^{co}Qcoh}_G({\rm V}(\mathcal{E}(\chi)),\chi,W+Q_s)$ to ${\rm Dcoh}_G(Z,\chi,W)$, we obtain an exact functor
$$i_*p^*:{\rm Dcoh}_G(Z,\chi,W)\rightarrow{\rm Dcoh}_G({\rm V}(\mathcal{E}(\chi)),\chi,W+Q_s).$$ 

Shipman proved that the above functor $i_*p^*$ is an equivalence when $G=\mathbb{G}_m$ trivially acts on $X$ and $W=0$ (see also \cite{isik}):

\begin{thm}[\cite{shipman} Theorem 3.4]\label{shipman}
The composition 
$$i_*p^*:{\rm Dcoh}_{\mathbb{G}_m}(Z,\chi_1,0)\xrightarrow{\sim}{\rm Dcoh}_{\mathbb{G}_m}({\rm V}(\mathcal{E}(\chi_1)),\chi_1,Q_s)$$ 
is an equivalence.
\end{thm}

The goal of this section is to show the following main result which is an analogy of the above theorem.

\begin{thm}\label{main result} Assume that $W|_Z:Z\rightarrow \mathbb{A}^1$ is flat.
The functor $$i_*p^*:{\rm Dcoh}_G(Z,\chi,W)\rightarrow{\rm Dcoh}_G({\rm V}(\mathcal{E}(\chi)),\chi,W+Q_s)$$
is an equivalence.
\end{thm}

\begin{rem}
Let $S$ be a smooth quasi-projective variety, and let $G$ be an affine reductive group acting on $S$.
Let $W:S\rightarrow \mathbb{A}^1$ be a $\chi:=\chi_1+\chi_2$-semi invariant non-constant regular function for some characters $\chi_i:G\rightarrow\mathbb{G}_m$. Let $X:={\rm V}(\mathcal{O}(\chi_1))\cong S\times \mathbb{A}^1_{x_1}$ be the $G$-vector bundle over $S$, and let $s\in\Gamma(X,\mathcal{O}(\chi_1))^G$ be the section corresponding to the $\chi_1$-semi invariant function $S\times \mathbb{A}^1_{x_1}\rightarrow\mathbb{A}^1$ which is defined as the projection $(s,x_1)\mapsto x_1$. Then, $S$ is isomorphic to the zero scheme of $s$, and the $G$-vector bundle ${\rm V}(\mathcal{O}(-\chi_1)(\chi))$ over $X$ is isomorphic to the $G$-variety $S\times\mathbb{A}^2_{x_1,x_2}$, where the $G$-weights of $x_i$  is given by $\chi_i$. By Theorem \ref{main result}, we have the following equivalence
$${\rm Dcoh}_G(S,\chi,W)\simeq{\rm Dcoh}_G(S\times\mathbb{A}^2_{x_1,x_2},\chi,W+x_1x_2).$$
  This kind of equivalence is know as {\it Kn\"orrer periodicity}, so the above theorem is considered as a generalization of  the original Kn\"orrer periodicity \cite[Theorem 3.1]{knorrer}.

\end{rem}

\subsection{Lemmas for the main theorem}
In this section, we provide some lemmas for the main result. Throughout this section, we consider the case when $G$ is trivial.

Set $$\omega_j:=\bigwedge^{r}(\mathcal{I}_Z/\mathcal{I}_Z^2)^{\vee}
\hspace{5mm}{\rm and}\hspace{5mm} \omega_i:=p^*\omega_j,$$ where $\mathcal{I}_Z$ is the ideal sheaf of $Z$ in $X$. These are  invertible sheaves on $Z$ and ${\rm V}(\mathcal{E})|_Z$ respectively. 
We define an exact functor $$i^!:{\rm D^{co}Qcoh}({\rm V}(\mathcal{E}),W+Q_s)\rightarrow{\rm D^{co}Qcoh}({\rm V}(\mathcal{E})|_Z,W)$$
as $i^!(-):={\bf L}i^*(-)\otimes\omega_i[-r]$. By \cite[Theorem 3.8]{efi-posi}, the above functor $i^!$ is right adjoint to $i_*:{\rm D^{co}Qcoh}({\rm V}(\mathcal{E})|_Z,W)\rightarrow{\rm D^{co}Qcoh}({\rm V}(\mathcal{E}),W+Q_s)$.
Let $$K:=K(q^*s,t)\in {\rm lfr}({\rm V}(\mathcal{E}), Q_s)$$ be the Koszul factorization of $q^*s\in \Gamma({\rm V}(\mathcal{E}), q^*\mathcal{E}^{\vee})$ and $t\in \Gamma({\rm V}(\mathcal{E}), q^*\mathcal{E})$, where $t$ is the tautological section.
 By abuse of notation, we denote by  $\mathcal{O}_Z$ the object in ${\rm coh}(Z,0)$ of the following form
 $$\Bigl(0\rightarrow \mathcal{O}_Z \rightarrow 0\Bigr).$$

\begin{lem}\label{koszul duality}
 Consider the case when $W=0$. We have isomorphisms $$i_*p^*(\mathcal{O}_Z)\cong K\hspace{5mm}{\rm and}\hspace{5mm} p_*i^!(K)\cong\mathcal{O}_Z$$ in ${\rm Dcoh}({\rm V}(\mathcal{E}),Q_s)$ and in ${\rm D^{co}Qcoh}(Z,0)$ respectively.
\begin{proof}
These isomorphisms follow from Lemma \ref{koszul lemma}. In particular, the former isomorphism is an immediate consequence. Note that $\omega_i\cong i^*\bigwedge^rq^*\mathcal{E}^{\vee}$. We obtain the latter isomorphism as follows;
$$p_*i^!(K)\cong p_*\textbf{L}i^*(\mathcal{O}_Z\otimes \bigwedge^rq^*\mathcal{E}^{\vee}[-r])\cong p_*\textbf{L}i^*(K^{\vee}) \cong  p_*\textbf{L}i^*(\mathcal{O}_{Z_{t^{\vee}}})\cong \mathcal{O}_Z,$$ 
where the last isomorphism follows from the fact that the zero section $Z\subset {\rm V}(\mathcal{E})$ is isomorphic to the fiber product of closed subschemes ${\rm V}(\mathcal{E})|_Z\hookrightarrow {\rm V}(\mathcal{E})$ and $Z_{t^{\vee}}\hookrightarrow{\rm V}(\mathcal{E})$.
\end{proof}
\end{lem}

\begin{lem}\label{main lemma} 
The functor $$i_*p^*:{\rm D^{co}Qcoh}(Z,W)\rightarrow{\rm D^{co}Qcoh}({\rm V}(\mathcal{E}),W+Q_s)$$
is  fully faithful.

\begin{proof}
The functors $i_*p^*$ and $p_*i^!$ can be represented as integral functors
$$i_*p^*\cong\Phi_{k_*\mathcal{O}_{{\rm V}(\mathcal{E})|_Z}}\hspace{3mm}{\rm and}\hspace{3mm}p_*i^!\cong\Phi_{k_*\omega_i[-r]},$$
where $k:=p\times i:{\rm V}(\mathcal{E})|_Z\rightarrow Z\times {\rm V}(\mathcal{E})$ and kernels $\mathcal{O}_{{\rm V}(\mathcal{E})|_Z}$ and $\omega_i[-r]$ are objects in ${\rm Dcoh}({\rm V}(\mathcal{E})|_Z,0)$. By easy computation, we see that there exists an object $P\in {\rm D^{co}Qcoh}(Z,0)$ such that $p_*i^!\circ i_*p^*\cong\Phi_{\Delta_*P}\cong(-)\otimes P$, where $\Delta:Z\rightarrow Z\times Z$ is the diagonal embedding. Substituting  $W=0$, by Lemma \ref{koszul duality}, we have an isomorphism $P\cong\mathcal{O}_Z$. But $P$ doesn't depend on the function $W$. Hence, for any $W$, we have an isomorphism of functors $p_*i^!\circ i_*p^*\cong\Phi_{\Delta_*P}\cong{\rm id}_{{\rm D^{co}Qcoh}(Z,W)}$. By the following lemma, this implies that the functor $i_*p^*:{\rm D^{co}Qcoh}(Z,W)\rightarrow{\rm D^{co}Qcoh}({\rm V}(\mathcal{E}),W+Q_s)$ is fully faithful.
\end{proof}
\end{lem}

The following lemma is an opposite version of \cite[Lemma 1.1.1]{johnstone}.\footnote{I would like to thank Timothy Logvinenko for informing me about this lemma.} We give a proof for the reader's convenience.


\begin{lem}[cf. \cite{johnstone} Lemma 1.1.1]
Let  $F:\mathcal{A}\rightarrow\mathcal{B}$ be a functor between (arbitrary) categories, and suppose that $F$ admits a right adjoint functor $G:\mathcal{B}\rightarrow\mathcal{A}$. Then, if there exists an isomorphism of functors $\alpha:{\rm id}_{\mathcal{A}}\xrightarrow{\sim}GF$ ($\alpha$ is not necessarily the adjunction morphism), then $F$ is fully faithful.
\begin{proof}
The isomorphism  $\alpha$ implies that the following composition of maps is bijective;
$${\rm Hom}(A,A')\xrightarrow{F}{\rm Hom}(F(A),F(A'))\xrightarrow{G}{\rm Hom}(GF(A),GF(A')).$$
Hence it is enough to show that $G$ is fully faithful on the image of $F$. Since the above composition is bijective, $G$ is full on the image of $F$.  Let $\varepsilon:{\rm id}_{\mathcal{A}}\rightarrow GF$ and $\delta: FG\rightarrow{\rm id}_{\mathcal{B}}$ be the adjunction morphisms. For any $f\in {\rm Hom}(F(A),F(A'))$ we have 
$$\delta_{F(A')}\circ FG(f)\circ F(\varepsilon_A)=f\circ \delta_{F(A)}\circ F(\varepsilon_A)=f,$$
where the first equation follows from the functoriality of $\delta$ and the the second equation follows from the property of the adjunction morphisms.
Hence the following  diagram is commutative
\[\xymatrix{
{\rm Hom}(F(A),F(A'))\ar@{=}[rr]\ar[d]_{G}&&{\rm Hom}(F(A),F(A'))\\
{\rm Hom}(GF(A),GF(A'))\ar[rr]^{F}&&{\rm Hom}(FGF(A),FGF(A'))\ar[u]_{\delta_{F(A')} \circ (-)\circ F(\varepsilon_A)},
}\]
and hence  $G$ is faithful on the image of $F$.
\end{proof}
\end{lem}


\subsection{Proof of the main theorem}
In this section, we prove the main theorem. Recall that $G$ is a reductive affine algebraic group acting on a smooth quasi-projective variety $X$. Since $X$ is smooth, there is a $G$-equivariant ample line bundle on $X$. In what follows, we assume that $W|_Z:Z\rightarrow \mathbb{A}^1$ is flat.

At first, we consider relative singularity categories. Let $Z_0$, $V|_{Z_0}$ and $V_0$ be the fibers of $W:Z\rightarrow \mathbb{A}^1$, $W:{\rm V}(\mathcal{E}(\chi))|_Z\rightarrow \mathbb{A}^1$ and $W+Q_s:{\rm V}(\mathcal{E}(\chi))\rightarrow \mathbb{A}^1$ over $0\in \mathbb{A}^1$ respectively. Denote by $p_0:V|_{Z_0}\rightarrow Z_0$ and $i_0:V|_{Z_0}\rightarrow V_0$ the restrictions of $p$ and $i$ respectively. By \cite[Corollary 2.27]{kuz}, the following cartesian squares are exact
\[\xymatrix{
V|_{Z_0}\ar[rr]^{p_0}\ar[d]&&Z_0\ar[d]&&V|_{Z_0}\ar[rr]^{i_0}\ar[d]&&V_0\ar[d]\\
{\rm V}(\mathcal{E})|_Z\ar[rr]^{p}&&Z&&{\rm V}(\mathcal{E})|_Z\ar[rr]^{i}&&\rm{V}(\mathcal{E}).
}\]
Since $p$ and $i$ have finite flat dimensions,  we have exact functors of relative singularity categories
$${p_0}^{\circ}:{\rm D}^{\rm cosg}_G(Z_0/Z)\rightarrow {\rm D}^{\rm cosg}_G(V|_{Z_0}/{\rm V}(\mathcal{E})|_Z)$$
$${i_0}_{\circ}:{\rm D}^{\rm cosg}_G(V|_{Z_0}/{\rm V}(\mathcal{E})|_Z)\rightarrow {\rm D}^{\rm cosg}(V_0/{\rm V}(\mathcal{E}))={\rm D}^{\rm cosg}_G(V_0).$$
Then the following diagram is commutative
\[\xymatrix{
{\rm D}^{\rm cosg}_G(Z_0/Z)\ar[rr]^{{i_0}_{\circ}{p_0}^{\circ}}\ar[d]_{\Upsilon}&&{\rm D}^{\rm cosg}_G(V_0) \ar[d]^{\Upsilon}\\
{\rm D^{co}Qcoh}_G(Z,W)\ar[rr]^-{i_*p^*}&&{\rm D^{co}Qcoh}_G({\rm V}(\mathcal{E}(\chi)),W+Q_s).
}\]

Furthermore, we compactify $V_0$ and $V|_{Z_0}$. The compactifying technique appeared  in \cite{orlov1}.
Let $$P:=\mathbb{P}(\mathcal{E}(\chi)\oplus \mathcal{O}_{X})=\underline{\rm Proj}({\rm Sym(\mathcal{E}(\chi)\oplus \mathcal{O}_{X})^{\vee}})$$ be the projective space bundle over $X$ with a $G$-action induced by the equivariant structure of $\mathcal{E}(\chi)\oplus \mathcal{O}_{X}$. Then we have a natural equivariant open immersion 
$$l:{\rm V}(\mathcal{E}(\chi))\rightarrow P.$$
Denote by $l|_Z: {\rm V}(\mathcal{E}(\chi))|_Z\rightarrow P|_Z$ the pull-back of $l$ by the closed immersion $\overline{i}:P|_Z\rightarrow P$. Now we have the following cartesian square
\[\xymatrix{
{\rm V}(\mathcal{E}(\chi))|_Z\ar[rr]^{l|_Z}\ar[d]_{i}&&P|_Z\ar[d]^{\overline{i}}\\
\rm{V}(\mathcal{E}(\chi))\ar[rr]^{l}&&P.
}\]
Denote by $\overline{q}:P\rightarrow X$ the natural projection, and let $\overline{p}:=\overline{q}|_Z:P|_Z\rightarrow Z$ be the pull-back of $\overline{q}$ by the closed immersion $j:Z\rightarrow X$.
Let $P_0$ be the $G$-invariant subscheme of $P$ defined by the $G$-invariant section $\overline{s\oplus W}\in \Gamma(P,\mathcal{O}(1)(\chi))^G$ which is corresponding to the composition
$$\mathcal{O}_P\xrightarrow{\overline{q}^*(s\oplus W)}\overline{q}^*({\mathcal{E}\oplus\mathcal{O}}(\chi^{-1}))^{\vee}\xrightarrow{\sigma} \mathcal{O}_{P}(1)(\chi),$$
where  $\sigma$ is the canonical surjection, and let $P|_{Z_0}$ be the zero scheme defined by the invariant section $\overline{i}^*(\overline{s\oplus W})\in \Gamma(P|_Z,\mathcal{O}(1)(\chi))^G$.  Since the pull-back  of $\overline{s\oplus W}$ (resp. $\overline{i}^*(\overline{s\oplus W})$) by the open immersion $l$ (resp. $l|_Z$) is equal to $W+Q_s$ (resp. $W$), we have the following exact cartesian square
\[\xymatrix{
V|_{Z_0}\ar[rr]^{l|_{Z_0}}\ar[d]_{i_0}&&P|_{Z_0}\ar[d]^{\overline{i_0}}\\
V_0\ar[rr]^{l_0}&&P_0.
}\]
Denote by $\overline{p_0}:P|_{Z_0}\rightarrow Z_0$ be the pull-back of $\overline{p}:P|_Z\rightarrow Z$ by the closed immersion  $Z_0\rightarrow Z$. Since the morphisms $\overline{i_0}:P|_{Z_0}\rightarrow P_0$ and $\overline{p_0}:P|_{Z_0}\rightarrow Z_0$ have finite Tor dimensions, the direct images $\textbf{R}\overline{i_0}_*:{\rm D^b}({\rm coh}P|_{Z_0})\rightarrow {\rm D^b}({\rm coh}P_0)$ and $\textbf{R}\overline{p_0}_*:{\rm D^b}({\rm coh}P|_{Z_0})\rightarrow {\rm D^b}({\rm coh}Z_0)$ induce the following exact functors (cf. \cite[Proposition 2.7]{tt}),
$$\overline{i_0}_{\circ}:{\rm D}^{\rm sg}_G(P|_{Z_0})\rightarrow {\rm D}^{\rm sg}_G(P_0)$$
$$\overline{p_0}_{\circ}:{\rm D}^{\rm sg}_G(P|_{Z_0})\rightarrow {\rm D}^{\rm sg}_G(Z_0).$$
Now we have the following commutative diagram

\[\xymatrix{
{\rm D}^{\rm sg}_G(Z_0)\ar[rr]^{\overline{i_0}_{\circ}{\overline{p_0}}^{\circ}}\ar[d]_{\pi}&&{\rm D}^{\rm sg}_G(P_0)\ar[d]^{{l_0}^{\circ}}\\
{\rm D}^{\rm sg}_G(Z_0/Z)\ar[rr]^{{i_0}_{\circ}{p_0}^{\circ}}&&{\rm D}^{\rm sg}_G(V_0),
}\]
where the vertical arrow on the left side is a Verdier localization  by Proposition \ref{verdier sing}.

\begin{rem}
If $Z$ is smooth, the above vertical arrows are equivalences. Indeed,  in that case, the singular locus ${\rm Sing}(P_0)$ is contained in $V_0$, whence ${l_0}^{\circ}$ is an equivalence by a similar argument in the proof of \cite[Proposition 1.14]{orlov2}. The equivalence of $\pi$ follows from Remark \ref{rel equal usual}.
\end{rem}

Let $\overline{i_0}^!:{\rm D^b}({\rm coh}_G P_0)\rightarrow {\rm D^b}({\rm coh}_G P|_{Z_0})$ be the functor defined by  $$\overline{i_0}^!:=\textbf{L}\overline{i_0}^*(-)\otimes \bigwedge^r(\mathcal{I}/\mathcal{I}^2)^{\vee}[-r],$$
where $\mathcal{I}$ is the  ideal sheaf of $\overline{i_0}:P|_{Z_0}\hookrightarrow P_0$.
The functor $\overline{i_0}^!$ is a right adjoint functor of $\overline{i_0}_*:{\rm D^b}({\rm coh}_{G}P|_{Z_0})\rightarrow{\rm D^b}({\rm coh}_{G}P_0)$. Indeed, these functors are adjoint when $G$ is trivial by \cite[I\hspace{-.1em}I\hspace{-.1em}I Theorem 6.7, Corollary 7.3]{hartshorne}, and the isomorphism $${\rm Hom}(\overline{i_0}_*(A),B)\cong{\rm Hom}(A,\overline{i_0}^!(B)),$$
where $A\in{\rm D^b}({\rm coh}P|_{Z_0})$ and $B\in{\rm D^b}({\rm coh}P_0)$, commutes with $G$-actions on each vector space of morphisms by the property in \cite[I\hspace{-.1em}I\hspace{-.1em}I Proposition 6.9.c]{hartshorne}.  Hence  we see that $\overline{i_0}^!$ is right adjoint to  $\overline{i_0}_*$ by \cite[Lemma 2.2.8]{vgit}. Denote by $$\overline{i_0}^{\flat}:{\rm D}^{\rm sg}_G(P_0)\rightarrow {\rm D}^{\rm sg}_G(P|_{Z_0})$$ the functor induced by $\overline{i_0}^!$. By the above argument, we have the following  adjoint pair $$\overline{i_0}_{\circ}\dashv\overline{i_0}^{\flat}.$$
Similarly, we have a right adjoint functor $${i_0}^{\flat}:{\rm D}^{\rm cosg}_G(V_0)\rightarrow{\rm D}^{\rm cosg}_G(V|_{Z_0}/{\rm V}(\mathcal{E}(\chi))|_Z)$$ of ${i_0}_{\circ}:{\rm D}^{\rm cosg}_G(V|_{Z_0}/{\rm V}(\mathcal{E}(\chi))|_Z)\rightarrow{\rm D}^{\rm cosg}_G(V_0)$.\\

{\it Proof of Theorem \ref{main result}:}
We have the following commutative diagram
\[\xymatrix{
{\rm D}^{\rm sg}_G(Z_0/Z)\ar[rr]^{{i_0}_{\circ}{p_0}^{\circ}}\ar[d]_{\Upsilon}&&{\rm D}^{\rm sg}_G(V_0) \ar[d]^{\Upsilon}\\
{\rm Dcoh}_G(Z,\chi,W)\ar[rr]^-{i_*p^*}&&{\rm Dcoh}_G({\rm V}(\mathcal{E}(\chi)),\chi,W+Q_s),
}\]
where the vertical arrows are equivalences by Theorem \ref{relative sing}. Hence it suffices to show that the functor ${i_0}_{\circ}{p_0}^{\circ}:{\rm D}^{\rm sg}_G(Z_0/Z)\rightarrow{\rm D}^{\rm sg}_G(V_0)$ is an equivalence.

At first, we prove that the functor ${i_0}_{\circ}{p_0}^{\circ}:{\rm D}^{\rm cosg}_G(Z_0/Z)\rightarrow{\rm D}^{\rm cosg}_G(V_0)$ is fully faithful. Let $$\varepsilon_G:{\rm id}_{{\rm D}^{\rm cosg}_G(Z_0/Z)}\rightarrow {p_0}_{\circ}{i_0}^{\flat}\circ {i_0}_{\circ}{p_0}^{\circ}$$
be the adjunction morphism of the adjoint pair ${i_0}_{\circ}{p_0}^{\circ}\dashv{p_0}_{\circ}{i_0}^{\flat}$. It is enough to show that for any object $A\in {\rm D}^{\rm cosg}_G(Z_0/Z)$, the cone $C_G(A)$ of the morphism $\varepsilon_G(A):A\rightarrow {p_0}_{\circ}{i_0}^{\flat}\circ {i_0}_{\circ}{p_0}^{\circ}(A)$ is the zero object. But the object ${\rm Res}_G(C_G(A))$ is isomorphic to the cone $C(A)$ of the adjunction morphism of $\varepsilon({\rm Res}_G(A)):{\rm Res}_G(A)\rightarrow {p_0}_{\circ}{i_0}^{\flat}\circ {i_0}_{\circ}{p_0}^{\circ}({\rm Res}_G(A))$ of the adjoint pair of  functors between ${\rm D}^{\rm cosg}(Z_0/Z)$ and ${\rm D}^{\rm cosg}(V_0)$. Since we have the following commutative diagram
\[\xymatrix{
{\rm D}^{\rm cosg}(Z_0/Z)\ar[rr]^{{i_0}_{\circ}{p_0}^{\circ}}\ar[d]_{\Upsilon}&&{\rm D}^{\rm cosg}(V_0) \ar[d]^{\Upsilon}\\
{\rm D^{co}Qcoh}(Z,W)\ar[rr]^-{i_*p^*}&&{\rm D^{co}Qcoh}({\rm V}(\mathcal{E}),W+Q_s),
}\]
where the vertical arrows are equivalences by Theorem \ref{relative sing}, the functor ${i_0}_{\circ}{p_0}^{\circ}$ is fully faithful by Lemma \ref{main lemma}. This implies that  the object $C(A)$ is the zero object. Hence $C_G(A)$ is also the zero object since the restriction functor ${\rm Res}_G$ is faithful by Lemma \ref{faithful lemma}. Hence ${i_0}_{\circ}{p_0}^{\circ}:{\rm D}^{\rm cosg}_G(Z_0/Z)\rightarrow{\rm D}^{\rm cosg}_G(V_0)$ is fully faithful. This implies that ${i_0}_{\circ}{p_0}^{\circ}:{\rm D}^{\rm sg}_G(Z_0/Z)\rightarrow{\rm D}^{\rm sg}_G(V_0)$ is also fully faithful, since the natural inclusions ${\rm D}^{\rm sg}_G(Z_0/Z)\rightarrow{\rm D}^{\rm cosg}_G(Z_0/Z)$ and ${\rm D}^{\rm sg}_G(V_0)\rightarrow{\rm D}^{\rm cosg}_G(V_0)$ are fully faithful by Theorem \ref{relative sing} and Proposition \ref{support inclusion} (1).

It only remains to show that the functor ${i_0}_{\circ}{p_0}^{\circ}:{\rm D}^{\rm sg}_G(Z_0/Z)\rightarrow{\rm D}^{\rm sg}_G(V_0)$ is essentially surjective. Consider the following commutative diagram:
\[\xymatrix{
{\rm D}^{\rm sg}_G(Z_0)\ar[rr]^{\overline{i_0}_{\circ}\overline{p_0}^{\circ}}\ar[d]_{\pi}&&{\rm D}^{\rm sg}_G(P_0) \ar[d]^{ {l_0}^{\circ}}\\
{\rm D}^{\rm sg}_G(Z_0/Z)\ar[rr]^-{{i_0}_{\circ}{p_0}^{\circ}}&&{\rm D}^{\rm sg}_G(V_0)
}\]
By a similar argument as in the proof of  \cite[Lemma 1.11]{orlov2}, we see that  every object in ${\rm D}^{\rm sg}_G(V_0)$ is isomorphic to an object $F[k]$ for some $G$-equivariant coherent sheaf $F$ and for some integer $k\in\mathbb{Z}$. Hence the vertical arrow on the right hand side in the above diagram is essentially surjective, since for every object $E$ in ${\rm coh}_GV_0$ there exists an object $\overline{E}$ in ${\rm coh}_GP_0$ such that ${l_0}^*(\overline{E})\cong E$.  Thus, we only need to prove that $\overline{i_0}_{\circ}\overline{p_0}^{\circ}:{\rm D}^{\rm sg}_G(Z_0)\rightarrow{\rm D}^{\rm sg}_G(P_0)$ is essentially surjective.
To prove that, it is enough to show  that the right adjoint functor $\overline{p_0}_{\circ}\overline{i_0}^{\flat}:{\rm D}^{\rm sg}_G(P_0)\rightarrow{\rm D}^{\rm sg}_G(Z_0)$ is fully faithful. Since the restriction functor ${\rm Res}_G:{\rm D}^{\rm sg}_G(P_0)\rightarrow{\rm D}^{\rm sg}(P_0)$ is faithful by Lemma \ref{faithful lemma} and \cite[Proposition 3.8]{pv}, it follows from \cite[Theorem 2.1]{orlov1} that the adjunction $\overline{i_0}_{\circ}\overline{p_0}^{\circ}\circ\overline{p_0}_{\circ}\overline{i_0}^{\flat}\rightarrow{\rm id}_{{\rm D}^{\rm sg}_G(P_0)}$ is an isomorphism of functors  by a similar argument as in the proof of  the fully faithfulness of ${i_0}_{\circ}{p_0}^{\circ}:{\rm D}^{\rm cosg}_G(Z_0/Z)\rightarrow{\rm D}^{\rm cosg}_G(V_0)$ in the previous paragraph.\qed

\subsection{Cases when $W=0$.}

In the previous section, we prove the main result assuming that $W|_Z:Z\rightarrow \mathbb{A}^1$ is flat. In this section, we consider the cases when $W=0$. In this cases, using results in \cite{shipman}, we can show the following:

\vspace{1mm}
With notation as above, consider $\mathbb{G}_m\times G$-action on $X$ induced by the projection $\mathbb{G}_m\times G\rightarrow G$. Let $\theta:\mathbb{G}_m\times G\rightarrow \mathbb{G}_m$ be the character defined as the projection. Since the first factor of $\mathbb{G}_m\times G$ trivially acts on $X$, the $G$-equivariant locally free sheaf $\mathcal{E}$ has a natural  $\mathbb{G}_m\times G$-equivariant structure. 

\begin{prop} 
We have an equivalence
$${\rm D^b}({\rm coh}_GZ)\xrightarrow{\sim}{\rm Dcoh}_{\mathbb{G}_m\times G}({\rm V}(\mathcal{E}(\theta)),\theta,Q_s).$$

\begin{proof}
By a similar argument as in section 2.2, we obtain an equivalence $${\rm D^b}({\rm coh}_GZ)\xrightarrow{\sim}{\rm Dcoh}_{\mathbb{G}_m\times G}(Z,\theta,0).$$
Hence it is enough to show the functor $$i_*p^*:{\rm Dcoh}_{\mathbb{G}_m\times G}(Z,\theta,0)\rightarrow{\rm Dcoh}_{\mathbb{G}_m\times G}({\rm V}(\mathcal{E}(\theta)),\theta,Q_s)$$
is an equivalence.

By Lemma \ref{main lemma}, it follows that  $$i_*p^*:{\rm D^{co}Qcoh}_{\mathbb{G}_m}(Z,\chi_1,0)\rightarrow{\rm D^{co}Qcoh}_{\mathbb{G}_m}({\rm V}(\mathcal{E}(\chi_1)),\chi_1,Q_s)$$ is fully faithful since the forgetful functor ${\rm D^{co}Qcoh}_{\mathbb{G}_m}(Z,\chi_1,0)\rightarrow{\rm D^{co}Qcoh}(Z,0)$ is faithful. Furthermore, the above functor $i_*p^*$ is an equivalence since the right orthogonal of the image of the restricted functor $i_*p^*:{\rm Dcoh}_{\mathbb{G}_m}(Z,\chi_1,0)\rightarrow{\rm D^{co}Qcoh}_{\mathbb{G}_m}({\rm V}(\mathcal{E}(\chi_1)),\chi_1,Q_s)$ vanishes by the argument in \cite[Theorem 3.4]{shipman}. In particular, the right adjoint functor 
$$p_*i^!:{\rm D^{co}Qcoh}_{\mathbb{G}_m}({\rm V}(\mathcal{E}(\chi_1)),\chi_1,Q_s)\rightarrow{\rm D^{co}Qcoh}_{\mathbb{G}_m}(Z,\chi_1,0)$$ of $i_*p^*$ is also fully faithful. 

Next we will show that the functor
$$i_*p^*:{\rm D^{co}Qcoh}_{\mathbb{G}_m\times G}(Z,\theta,0)\rightarrow{\rm D^{co}Qcoh}_{\mathbb{G}_m\times G}({\rm V}(\mathcal{E}(\theta)),\theta,Q_s)$$
is an equivalence.
Let $$\varepsilon_{\mathbb{G}_m\times G}:{\rm id}_{{\rm D^{co}Qcoh}_{\mathbb{G}_m\times G}(Z,\theta,0)}\rightarrow p_*i^!\circ i_*p^*$$ be the adjunction morphism. To show that the functor $i_*p^*:{\rm D^{co}Qcoh}_{\mathbb{G}_m\times G}(Z,\theta,0)\rightarrow{\rm D^{co}Qcoh}_{\mathbb{G}_m\times G}({\rm V}(\mathcal{E}(\theta)),\theta,Q_s)$ is fully faithful, we will prove that the adjunction morphism $\varepsilon_{\mathbb{G}_m\times G}$ is an isomorphism of functors. For this, it suffices to show that for any object $F\in{\rm D^{co}Qcoh}_{\mathbb{G}_m\times G}(Z,\theta,0)$ the cone $C_{\mathbb{G}_m\times G}(F)$ of the morphism $\varepsilon_{\mathbb{G}_m\times G}(F):F\rightarrow p_*i^!\circ i_*p^*(F)$ is the zero object.
Recall that the categories ${\rm Qcoh}_{\mathbb{G}_m}Z$ and ${\rm Qcoh}_{\mathbb{G}_m\times G}Z$ are equivalent to the categories ${\rm Qcoh}[Z/\mathbb{G}_m]$ and ${\rm Qcoh}_G[Z/\mathbb{G}_m]$ respectively, where   $[Z/\mathbb{G}_m]$ denotes the quotient stack,  and we can consider the restriction and the induction functors for algebraic stacks as in section 2.6. Let $\pi_G:{\rm Qcoh}_{\mathbb{G}_m\times G}Z\rightarrow{\rm Qcoh}_{\mathbb{G}_m}Z$ be the functor corresponding to the restriction functor ${\rm Res}_G:{\rm Qcoh}_G[Z/\mathbb{G}_m]\rightarrow{\rm Qcoh}[Z/\mathbb{G}_m]$ via the equivalences ${\rm Qcoh}_{\mathbb{G}_m}Z\cong{\rm Qcoh}[Z/\mathbb{G}_m]$ and ${\rm Qcoh}_{\mathbb{G}_m\times G}Z\cong {\rm Qcoh}_G[Z/\mathbb{G}_m]$.
Then $\pi_G$ naturally induces the following exact functor 
$$\pi_G:{\rm D^{co}Qcoh}_{\mathbb{G}_m\times G}(Z,\theta,0)\rightarrow {\rm D^{co}Qcoh}_{\mathbb{G}_m}(Z,\chi_1,0),$$  
and  $\pi_G$ has the right adjoint  functor $\sigma_G:{\rm D^{co}Qcoh}_{\mathbb{G}_m}(Z,\chi_1,0)\rightarrow{\rm D^{co}Qcoh}_{\mathbb{G}_m\times G}(Z,\theta,0)$ induced by the induction functor. Since the argument  in the proof of Lemma \ref{reductive} works for algebraic stacks, the adjunction morphism  ${\rm id}\rightarrow \sigma_G\circ\pi_G$ is a split mono. Hence $\pi_G$ is faithful.
The object $\pi_G(C_{\mathbb{G}_m\times G}(F))$ is isomorphic to the cone $C_{\mathbb{G}_m}(F)$ of the adjunction morphism $\varepsilon_{\mathbb{G}_m}(\pi_G(F)):\pi_G(F)\rightarrow p_*i^!\circ i_*p^*(\pi_G(F))$, and $C_{\mathbb{G}_m}(F)$ is the zero object since the functor $i_*p^*:{\rm D^{co}Qcoh}_{\mathbb{G}_m}(Z,\chi_1,0)\rightarrow{\rm D^{co}Qcoh}_{\mathbb{G}_m}({\rm V}(\mathcal{E}(\chi_1)),\chi_1,Q_s)$ is fully faithful.  Hence  we see that the object $C_{\mathbb{G}_m\times G}(F)$ is also the zero object since $\pi_G$ is faithful. 
By an identical argument, we see that the right adjoint functor $$p_*i^!:{\rm D^{co}Qcoh}_{\mathbb{G}_m\times G}({\rm V}(\mathcal{E}(\theta)),\theta,Q_s)\rightarrow{\rm D^{co}Qcoh}_{\mathbb{G}_m\times G}(Z,\theta,0)$$ is also fully faithful. Hence  the functor $$i_*p^*:{\rm D^{co}Qcoh}_{\mathbb{G}_m\times G}(Z,\theta,0)\rightarrow{\rm D^{co}Qcoh}_{\mathbb{G}_m\times G}({\rm V}(\mathcal{E}(\theta)),\theta,Q_s)$$
is an equivalence.

By Proposition \ref{support inclusion} (1),  we see that the equivalence $i_*p^*:{\rm D^{co}Qcoh}_{\mathbb{G}_m\times G}(Z,\theta,0)\rightarrow{\rm D^{co}Qcoh}_{\mathbb{G}_m\times G}({\rm V}(\mathcal{E}(\theta)),\theta,Q_s)$ induces an equivalence of  the compact objects $$i_*p^*:\overline{{\rm Dcoh}_{\mathbb{G}_m\times G}(Z,\theta,0)}\rightarrow\overline{{\rm Dcoh}_{\mathbb{G}_m\times G}({\rm V}(\mathcal{E}(\theta)),\theta,Q_s)},$$
where $\overline{(-)}$ denotes the idempotent completion of $(-)$. But ${\rm Dcoh}_{\mathbb{G}_m\times G}(Z,\theta,0)$ on the left hand side is already idempotent complete since it is equivalent to ${\rm D^b}({\rm coh}_GZ)$. Hence the functor 
$$i_*p^*:{\rm Dcoh}_{\mathbb{G}_m\times G}(Z,\theta,0)\rightarrow{\rm Dcoh}_{\mathbb{G}_m\times G}({\rm V}(\mathcal{E}(\theta)),\theta,Q_s)$$
is an equivalence.
\end{proof}
\end{prop}

\section{Orlov's theorem for gauged LG models}

In this section, we obtain a gauged LG version of the following  theorem of Orlov.

\begin{thm}[\cite{orlov3} Theorem 40]\label{orlov sod}
Let $X\subset \mathbb{P}^{N-1}_k$ be the hypersurface defined by a section $f\in \Gamma(\mathbb{P}^{N-1}_k,\mathcal{O}(d))$. Denote by $F$ the corresponding homogeneous polynomial.\\
$(1)$ If $d<N$, there is a semi-orthogonal decomposition 
$${\rm D^b}({\rm coh}X)=\langle \mathcal{O}_X(d-N+1),...,\mathcal{O}_X,{\rm Dcoh}_{\mathbb{G}_m}(\mathbb{A}^N_k,\chi_{d},F)\rangle.$$
$(2)$ If $d=N$, there is an equivalence
$${\rm D^b}({\rm coh}X)\cong {\rm Dcoh}_{\mathbb{G}_m}(\mathbb{A}^N_k,\chi_{d},F).$$
$(3)$ If $d>N$, there is a semi-orthogonal decomposition
$${\rm Dcoh}_{\mathbb{G}_m}(\mathbb{A}^N_k,\chi_{d},F)=\langle k,...,k(N-d+1),{\rm D^b}({\rm coh}X)\rangle.$$
\end{thm}

We combine the main result with the theory of variations of GIT quotients to obtain a gauged LG version of the above theorem. For the theory of variations of GIT quotients, see \cite{vgit} or \cite[Section 2]{degree d}. This kind of approach to Orlov's theorem appeared in \cite{shipman}, \cite{vgit}, and \cite{degree d}, and our argument is similar to the one in \cite[Section 3]{degree d}.

\vspace{4mm}
Let $S$ be a smooth quasi-projective variety with $\mathbb{G}_m$-action, and set
$$Q:=S\times \mathbb{A}^N\times\mathbb{A}^1.$$
 For $i=1, 2$, set $G_i:=\mathbb{G}_m$, and let $G:=G_1\times G_2$. For a positive integer $d>1$, we define a $G$-action on $Q$ as follows;
$$G\times Q\ni (g_1,g_2) \times(s,v_1,...v_N,u)\mapsto (g_2\cdot s,g_1v_1,...,g_1v_N,g_1^{-d}g_2 u)\in Q,$$
where the action $\cdot$ is  the original $\mathbb{G}_m$-action on $S$. Let $\lambda:\mathbb{G}_m\rightarrow G$ be the character defined by $\lambda(a):=(a,1)$. Denote by  $Z_{\lambda}$ the fixed locus of $\lambda$-action on $Q$. Then $Z_{\lambda}$ coincides with the zero section $S\times{0}\times{0}\subset Q$. Furthermore, set 
$S_+:=\{q\in Q\mid \displaystyle{\lim_{a\to 0}}\lambda(a)q\in Z_{\lambda}\}$ and 
$S_-:=\{q\in Q\mid  \displaystyle{\lim_{a\to 0}}\lambda(a)^{-1}q\in Z_{\lambda}\}.$
Then $$S_+=S\times \mathbb{A}^N\times0 \hspace{5mm}{\rm and }\hspace{5mm}S_-=S\times 0\times\mathbb{A}^1.$$
Denote by $Q_+$ (resp. $Q_-$) be the complement of $S_+$ (resp. $S_-$) in $Q$.
Then the stratifications
$$Q=Q_+\sqcup S_+\hspace{5mm}{\rm and }\hspace{5mm}Q=Q_-\sqcup S_-$$
are elementary wall crossings in the sense of \cite{vgit}.

Let $W: S\rightarrow \mathbb{A}^1$ be a $\chi_1$-semi invariant function which is flat. Let $f\in \Gamma(\mathbb{P}^{N-1}_S,\mathcal{O}(d))^{\mathbb{G}_m}$ be a non-zero $\mathbb{G}_m$-invariant section, and denote by $F:\mathbb{A}_S^N\rightarrow\mathbb{A}^1$ the corresponding regular function.
Since $Q$ is the trivial line bundle over $\mathbb{A}_S^N$, the function $F$ induces a  regular function $\widetilde{F}:Q\rightarrow\mathbb{A}^1$. 
 Then the  function $$W+\widetilde{F}:Q\rightarrow \mathbb{A}^1$$ is a $\chi_{0,1}$-semi invariant regular function, where $W$ is the pull-back of $W:S\rightarrow\mathbb{A}^1$ by the projection $Q\rightarrow S$, and $\chi_{0,1}:G\rightarrow\mathbb{G}_m$ is the character defined by $\chi_{0,1}(g_1,g_2):=g_2$. By  \cite[Lemma 3.4.4]{vgit} and \cite[Theorem 3.5.2]{vgit}, we have the following:

\begin{prop}\label{orlov prop}
Let $t_{\pm}$ be the $\lambda$-weight of the restriction of relative canonical bundle $\omega_{S_{\pm}/Q}$ to $Z_{\lambda}$, and set $\mu:=-t_{+}+t_{-}$. Let $\chi:G\rightarrow \mathbb{G}_m$ be the character defined by $\chi(g_1,g_2):=g_1g_2$. \\
$(1)$ If $\mu<0$, there exist  fully faithful functors $$\Upsilon_-:{\rm Dcoh}_{G/\lambda}(Z_{\lambda},\chi_1,W+\widetilde{F})\rightarrow {\rm Dcoh}_G(Q_-,\chi_{0,1},W+\widetilde{F})$$  $$\Phi_-:{\rm Dcoh}_G(Q_+,\chi_{0,1},W+\widetilde{F})\rightarrow{\rm Dcoh}_G(Q_-,\chi_{0,1},W+\widetilde{F}),$$ and we have the following semi-orthogonal decomposition
$${\rm Dcoh}_G(Q_-,\chi_{0,1},W+\widetilde{F})=\langle \Upsilon_-(\mu+1),...,\Upsilon_-,\Phi_-({\rm Dcoh}_G(Q_+,\chi_{0,1},W+\widetilde{F}))\rangle,$$
where we denote by $\Upsilon_-(n)$ the the essential image of the composition $(-)\otimes\mathcal{O}(\chi^n)\circ \Upsilon_-$.
$(2)$ If $\mu=0$, we have an equivalence
$${\rm Dcoh}_G(Q_-,\chi_{0,1},W+\widetilde{F})\cong{\rm Dcoh}_G(Q_+,\chi_{0,1},W+\widetilde{F}).$$
$(3)$ If $\mu>0$, there exist  fully faithful functors $$\Upsilon_+:{\rm Dcoh}_{G/\lambda}(Z_{\lambda},\chi_1,W+\widetilde{F})\rightarrow {\rm Dcoh}_G(Q_+,\chi_{0,1},W+\widetilde{F})$$  $$\Phi_+:{\rm Dcoh}_G(Q_-,\chi_{0,1},W+\widetilde{F})\rightarrow{\rm Dcoh}_G(Q_+,\chi_{0,1},W+\widetilde{F}),$$ and we have the following semi-orthogonal decomposition
$${\rm Dcoh}_G(Q_+,\chi_{0,1},W+\widetilde{F})=\langle \Upsilon_+,...,\Upsilon_+(-\mu+1),\Phi_+({\rm Dcoh}_G(Q_-,\chi_{0,1},W+\widetilde{F}))\rangle,$$
where we denote by $\Upsilon_+(n)$ the the essential image of the composition $(-)\otimes\mathcal{O}(\chi^n)\circ \Upsilon_+$.

\end{prop}

\vspace{5mm}
Since $Z_{\lambda}=S\times0\times0$, the function $\widetilde{F}$ vanishes on $Z_{\lambda}\subset Q$. Hence we have 
$${\rm Dcoh}_{G/\lambda}(Z_{\lambda},\chi_1,W+\widetilde{F})\cong{\rm Dcoh}_{\mathbb{G}_m}(S,\chi_1,W).$$

Next, we have $$Q_-=S\times\mathbb{A}^N\setminus0\times\mathbb{A}^1.$$
Since $F|_{S\times\mathbb{A}^N\setminus0}\in\Gamma(S\times\mathbb{A}^N\setminus0,\mathcal{O}(\chi_{-d,0})^{\vee})^G$ and $Q_-={\rm V}(\mathcal{O}(\chi_{-d,1}))$, Theorem \ref{main result} implies the following equivalence;
$${\rm Dcoh}_G(Q_-,\chi_{0,1},W+\widetilde{F})\cong{\rm Dcoh}_G(Z,\chi_{0,1},W),$$
where $Z\subset S\times\mathbb{A}^N\setminus0$ is the zero scheme of $F$.  Moreover, the quotient stack $[Z/G_1]$ is isomorphic to the hypersurface $X$ in the projective space bundle $ \mathbb{P}^{N-1}_S$ over $S$ defined by the invariant section $f\in\Gamma(\mathbb{P}^{N-1}_S,\mathcal{O}(d))^{G_2}$. Hence we have an equivalence
$${\rm Dcoh}_G(Z,\chi_{0,1},W)\cong{\rm Dcoh}_{G_2}(X,\chi_1,W).$$

On the other hand, we have 
$$Q_+=S\times\mathbb{A}^N\times\mathbb{A}^1\setminus0.$$
We consider another action of $G$ on $Q_+$ as follows;
$$G\times Q_+\ni(g_1,g_2)\times(s,v,u)\mapsto (g_1^d\cdot s,g_1v,g_1^{-d}g_2u)\in Q_+.$$
We denote by $\widetilde{Q_+}$ the new $G$-variety. Then we have a $G$-equivariant isomorphism
$$\varphi:\widetilde{Q_+}\xrightarrow{\sim} Q_+,$$
given by $\varphi(s,v,u):=(u\cdot s,v,u),$ where $u\in\mathbb{A}^1\setminus0$ is considered as a point in $\mathbb{G}_m$.
Since $G_2$ trivially acts on the first two components $S\times \mathbb{A}^N$ of $\widetilde{Q_+}$, we have 
$$[\widetilde{Q_+}/G_2]\cong S\times\mathbb{A}^N\times[\mathbb{A}^1\setminus0/G_2]\cong \mathbb{A}^N_S.$$
Hence we have the following equivalence
$${\rm Dcoh}_G(Q_+,\chi_{0,1},W+\widetilde{F})\cong{\rm Dcoh}_{G_1}(\mathbb{A}^N_S,\chi_d,W+F),$$
where, on the right hand side, $G_1$-action is given by the  following
$$G_1\times S\times \mathbb{A}^N\ni g_1\times (s,v)\mapsto (g_1^d\cdot s, g_1v).$$

Finally, note that $\mu=d-N$ and that the twisting by the $G$-equivariant invertible sheaf $\mathcal{O}(\chi)$  corresponds to the twisting, in ${\rm Dcoh}_{G_2}(X,\chi_1,W)$, by the $G_2$-equivariant  invertible sheaf $\mathcal{O}(1)$ on $X$ which is the pull-back of the tautological $G_2$-equivariant  invertible sheaf on $\mathbb{P}^{N-1}_S$. Combining Proposition \ref{orlov prop} and the above argument, we obtain the following gauged LG version of the Orlov's theorem:

\vspace{5mm}
Let $S$ be a smooth quasi-projective variety with a $\mathbb{G}_m$-action, and let $W:S\rightarrow \mathbb{A}^1$ be a  $\chi_1$-semi invariant regular function which is flat. Consider $\mathbb{G}_m$-actions on $\mathbb{A}^N_S$ and on $\mathbb{P}^{N-1}_S$ given by 
$$\mathbb{G}_m\times\mathbb{A}^N_S\ni t\times (s,v_1,...,v_N)\mapsto (t^d\cdot s, tv_1,...tv_N)\in\mathbb{A}_S^N$$
$$\mathbb{G}_m\times\mathbb{P}^{N-1}_S\ni t\times (s,v_1:...:v_N)\mapsto (t\cdot s, v_1:...:v_N)\in\mathbb{P}^{N-1}_S.$$ 

\begin{thm}
  For $d>1$, let $f\in \Gamma(\mathbb{P}^{N-1}_S,\mathcal{O}(d))^{\mathbb{G}_m}$ be a non-zero invariant  section, and let $F:\mathbb{A}_S^N\rightarrow\mathbb{A}^1$ be the corresponding $\chi_d$-semi invariant regular function. Let $X\subset\mathbb{P}^{N-1}_S$ be the hypersurface defined by $f$, and assume that the morphism $W|_X$ is flat.\\
 $(1)$ If $d<N$, there are fully faithful functors 
 $$\Phi:{\rm Dcoh}_{\mathbb{G}_m}(\mathbb{A}^N_S,\chi_d,W+F)\rightarrow {\rm Dcoh}_{\mathbb{G}_m}(X,\chi_1,W)$$
 $$\Upsilon:{\rm Dcoh}_{\mathbb{G}_m}(S,\chi_1,W)\rightarrow {\rm Dcoh}_{\mathbb{G}_m}(X,\chi_1,W),$$
 and there is a semi-orthogonal decomposition
$${\rm Dcoh}_{\mathbb{G}_m}(X,\chi_1,W)=\langle\Upsilon_{d-N+1},...,\Upsilon_0,\Phi( {\rm Dcoh}_{\mathbb{G}_m}(\mathbb{A}^N_S,\chi_d,W+F))\rangle,$$
where $\Upsilon_i$ denotes the essential image of the composition $(-)\otimes\mathcal{O}(i)\circ\Upsilon$.\\

\begin{flushleft}$(2)$ If $d=N$, we have an equivalence\end{flushleft}
$${\rm Dcoh}_{\mathbb{G}_m}(X,\chi_1,W)\cong{\rm Dcoh}_{\mathbb{G}_m}(\mathbb{A}^N_S,\chi_d,W+F).$$\\
$(3)$ If $d>N$, there are fully faithful functors 
 $$\Psi:{\rm Dcoh}_{\mathbb{G}_m}(X,\chi_1,W)\rightarrow {\rm Dcoh}_{\mathbb{G}_m}(\mathbb{A}^N_S,\chi_d,W+F)$$
 $$\Upsilon:{\rm Dcoh}_{\mathbb{G}_m}(S,\chi_1,W)\rightarrow {\rm Dcoh}_{\mathbb{G}_m}(\mathbb{A}^N_S,\chi_d,W+F),$$ and there is a semi-orthogonal decomposition
$${\rm Dcoh}_{\mathbb{G}_m}(\mathbb{A}^N_S,\chi_d,W+F)=\langle\Upsilon_0,...,\Upsilon_{N-d+1},\Psi({\rm Dcoh}_{\mathbb{G}_m}(X,\chi_1,W))\rangle,$$
 where $\Upsilon_i$ denotes the essential image of the composition $(-)\otimes\mathcal{O}(\chi_i)\circ\Upsilon$.
\end{thm}

\begin{rem}
(1) We can view  Orlov's Theorem \ref{orlov sod} as the case when $S={\rm Spec}\,k$ and $W=0$ in the above theorem.\\
(2) If $N>1$, the assumption that $W|_X$ is flat  is  satisfied whenever $W:S\rightarrow\mathbb{A}^1$ is flat.\\
(3) For positive integers $a_1,...,a_N$, applying the similar argument to the $G$-action on $Q$ defined by 
$$G\times Q\ni (g_1,g_2) \times(s,v_1,...v_N,u)\mapsto (g_2\cdot s,g_1^{a_1}v_1,...,g_1^{a_N}v_N,g_1^{-d}g_2 u)\in Q,$$
we can obtain the similar result for the hypersurface $X$ in weighted projective stack bundle $\mathbb{P}^{N-1}_S(a_1,...,a_{N}):=[S\times\mathbb{A}^N\setminus0/G_1]$ over $S$ defined by the section corresponding to a $G_1$-invariant section $F\in\Gamma(\mathbb{A}^N_S,\mathcal{O}(\chi_d))^{G_1}$.\\
(4) Of course, Orlov's theorem in \cite{orlov3} is much more general. It covers noncommutative situations unlike our setting.
\end{rem}

\vspace{4mm}
\address{Department of Mathematics and Information Sciences, Tokyo Metropolitan University, 1-1 Minamiohsawa, Hachioji-shi, Tokyo, 192-0397, Japan}\\
{\it E-mail address}: \email{yuki-hirano@ed.tmu.ac.jp}

\end{document}